# Random walks on binary strings applied to the somatic hypermutation of B-cells

Irene Balelli · Vuk Milišić · Gilles Wainrib

**Abstract** Within the germinal center in follicles, B-cells proliferate, mutate and differentiate, while being submitted to a powerful selection : a micro-evolutionary mechanism at the heart of adaptive immunity. A new foreign pathogen is confronted to our immune system, the mutation mechanism that allows B-cells to adapt to it is called *somatic hypermutation* : a programmed process of mutation affecting B-cell receptors at extremely high rate. By considering random walks on graphs, we introduce and analyze a simplified mathematical model in order to understand this extremely efficient learning process. The structure of the graph reflects the choice of the mutation rule. We focus on the impact of this choice on typical time-scales of the graphs' exploration. We derive explicit formulas to evaluate the expected hitting time to cover a given Hamming distance on the graphs under consideration. This characterizes the efficiency of these processes in driving antibody affinity maturation. In a further step we present a biologically more involved model and discuss its numerical outputs within our mathematical framework. We provide as well limitations and possible extensions of our approach.

I. Balelli · V. Milišić
Université Paris 13, Institut Galilée, Département de Mathématiques.
99, avenue Jean-Baptiste Clément 93430 - Villetaneuse - France.
Tel.: +33-1-49-40-36-39 ; +33-1-49-40-35-91
Fax: +33-1-49-40-35-68
E-mail: balelli@math.univ-paris13.fr
E-mail: milisic@math.univ-paris13.fr

G. Wainrib
Ecole Normale Supérieure, Département d'Informatique.
45 rue d'Ulm, 75005 - Paris - France.
E-mail: gilles.wainrib@ens.fr

**Contents**

# 1 Introduction

Understanding the role and functional implication of mutations is a central question in biological evolutionary theory [17,57,24,15], but also for the study of evolutionary algorithms [3]. Beyond the mutation rate, which is naturally an important parameter, our aim in this article is to highlight the role of various mutation rules on the exploration of the space of traits. In our mathematical framework, configurations are represented as vertices of a graph which are connected if there exists a mutation allowing to pass from one state to another. We are mainly interested in understanding the characteristic time-scales for the exploration of the state-space as a function of the mutation rule. To this end, we relate mutation rules with specific graph topologies and build upon random walks on graphs and spectral graph theories to analyze resulting time-scales.

More precisely, beyond general theoretical results, we are particularly interested to apply our framework to the B-cell affinity maturation in Germinal Centers (GCs). The adaptive immune system is able to create a specific response against almost any kind of pathogens penetrating our organism and inflict a disease. This task is performed by the production of high affinity antigen-specific antibodies. These proteins are produced by B-lymphocytes which are submitted to a learning process improving their affinity to recognize a particular antigen. This process is called Antibody Affinity Maturation (AAM) and takes place in GCs [40]. Even if substantial progress has been made in adaptative immunology, since somatic hypermutation was discovered by the nobel price Susumu Tonegawa [54] in 1987, there are still facts that remain unclear about the GC reaction and the exact dynamics of AAM. Indeed, it seems difficult to make exact measurements of the antigenic repertoire *in vivo*, or to have precise spatial and temporal data about lymphocytes during the GC reaction, or to understand the exact dynamic of mutation and selection of B-cells while they are submitted to AAM (*e.g.* [16,43]). Nevertheless, some refined techniques start to be available [21], showing possible correlations between proliferation and mutation rates with respect to B-cells' affinity. This provides further motivation for setting appropriate mathematical frameworks to describe such systems.

The affinity of a B-cell is biologically observed as a matching between the B-cell receptor (BCR) and the antigen. We aim at understanding how mu-tation rules allow to explore the state-space of all possible configurations of BCRs. The mutational mechanism that B-cells undergo in order to improve their affinity is called Somatic Hypermutation (SHM): it targets, at a very high rate, the DNA encoding for the specific portion of the BCR involved in the binding with the antigen, called Variable (V) regions. SHM can introduce mutations at all four nucleotides, and mutation hot-spots have been identified [53]. The effect of these mutations on the BCR, once expressed on the outer surface of B-cells, is very complex, as the substitution of a single amino-acid can modify the geometrical structure of the BCR, creating or deleting bonds (see [1], Chapter 4, for more details about the crystal structure of BCRs and their binding with antigens). All these features contribute in many different ways in the exploration of possible BCRs' configurations.

Although the mutation process occurs at the level of the DNA, the result of the mutation can be summarized by the modification of the amino-acid string composing the BCR. Therefore, in the present paper, we consider effective mutations directly at the level of the amino-acid string (Section 4.3). However, the structure of our mathematical model can be left substantially unchanged even in the hypothesis that we are considering mutations on the DNA, which leads to a modify the definition of affinity and the size of the state-space.

There already exists a certain number of mathematical models and results about GC reaction and AAM. In particular, [30,31] proposed deterministic population modeling of SHM and AAM, considering for instance the hypothesis of recycling mechanisms during GC's reaction, later investigated by experiments [55]. In [42,44], the authors introduced and discussed several immunological problems, such as the size of the repertoire, or the strength of antigen-antibody binding, or the pourcentage of recycling. They provide suitable mathematical tools, using both deterministic and probabilistic approaches, together with numerical simulations. More recently, biologically very detailed models of GCs were proposed [37], using, for instance, agent-based models [38], mostly analyzed through extensive numerical simulations. Our aim here is not to build a very detailed and sophisticated model, but rather to contribute to the theoretical foundation of adaptive immunity modeling through the mathematical analysis of generic mutation models on graphs. So far, this approach has not been developed and applied to GC reaction and AAM modeling. In particular, this framework enables the study of various mutation rules, such as for instance, affinity-dependent mutations, which are currently debated in the biological literature [21].

Beyond the fundamental understanding of physiological processes and their associated pathologies, this research is related to important biotechnological applications, such as the synthetic production of specific antibodies for drugs, vaccines or cancer immunotherapy [20,50], since this production process involves the selection of high affinity peptides and requires smart methods to generate an appropriate diversity, and also to the theoretical understanding of bio-inspired algorithms such as in [9].

In this article, we consider pure mutational models obtained as random walks on graphs given by alterations of the edge set of the $N$-dimensional hypercube. We are mostly interested in understanding the variation of hitting times as a function of the underlying graphs, hence relating mutation rules to the characteristic time-scales of the process. Section 2 contains mainly results which are based on random walks theory (*e.g.* [41,39,46]) and, more specifically, random walks on graphs (*e.g.* [36,2]). This is a topic of active research due to the great number of important applications in recent years, such as graph clustering [48], ranking algorithms for search-engines [6,27], or social network modeling [29,22,32].

For the sake of clarity, in Section 2 we start with the most basic mutational model which is the simple random walk on the $N$-dimensional hypercube (*e.g.* [13,25,12,56]). We set notations useful in order to define our models, then we briefly overview some basic properties of random walks on graphs, and establish particular results in the case of the hypercube. In Section 3, we study several mutation rules and their effects on the structure of the graph and, consequently, its associated random walk, in particular in terms of the hitting times. We use both spectral and probabilistic methods. We especially focus on two mutation rules that are the combination of simpler ones: the class switch of 1 or 2-length strings (Section 3.1.3), where the mutation rule depends on the distance to the target, and the mutation rule which allows to do more than a single mutation at each step (Section 3.1.4). Table 1 in Section 3.2 summarizes the main results of Section 2 and 3 : we display expected times to reach some position of the graph, as a function of each mutation rule. Finally, Section 4 is dedicated to modeling aspects and discussions about possible extensions and limitations of the proposed framework.

## 2 A basic mutational model

In this section we start by setting up the general mathematical framework, which we will keep to pattern and study all mutational mechanisms discussed in the current section and in Section 3. Hence, we state a basic mutational model. The choice of this environnement is motivated by the modeling of amino-acids chains and their modifications during SHM. It is for this reason that we often recall some biological facts and refer to BCRs and antigens to provide motivation. Despite this, this framework seems to us the simplest and most adaptable one to study different mutational rules in a more general evolutionary context.

We suppose it is possible to classify the amino-acids, which determine the chemical properties of both BCR and antigen, into 2 classes denoted by 0 and 1 respectively (they could represent amino-acids negatively and positively charged respectively). Henceforth BCRs and antigen are represented by binary strings of same fixed length $N$, hence, the state-space of all possible BCR con-figurations is $\{0,1\}^N$. We will give some more details about these hypotheses in Section 4.3.

**Definition 1** We denote by $\mathcal{H}_N$ the standard $N$-dimensional hypercube. BCR and antigen configurations are represented by vertices of $\mathcal{H}_N$, denoted by $\mathbf{x}_i$ with $1 \leq i \leq 2^N$, or sometimes simply by their indices. We denote the antigen target vertex by $\overline{\mathbf{x}}$: it is given at the beginning of the process and never changes.

We suppose that there is a single B-cell entering the GC reaction. The configuration of its receptor is denoted by $\mathbf{X}_0$. If $\mathbf{X}_t$ is the configuration of the BCR after $t$ mutations, then depending on the mutational rule, one or more bits in $\mathbf{X}_t$ can change after the next mutation. This gives rise to a Random Walk (RW) on $\{0,1\}^N$, where a mutation on the BCR corresponds to a jump to a neighbor node. Of course, the definition of neighbors changes depending on the mutation rules we introduce (we specify neighborhood each time we discuss a new mutation rule). In a general way:

**Definition 2** Given $\mathbf{x}_i, \mathbf{x}_j \in \{0,1\}^N$, we say that $\mathbf{x}_i$ and $\mathbf{x}_j$ are neighbors, and denote $\mathbf{x}_i \sim \mathbf{x}_j$, if there exists at least one edge (or loop) between them.

As far as the complementarity is concerned, we have to make a further simplification. As we have already discuss in the Introduction, the tridimensional structure of the BCR is hard to model. For this reason we consider a linear contact, *i.e.* positively charged amino-acids are complementary to negatively charged ones when they are at the same position within the binary string. For the sake of simplicity, we state that 0 matches with 0 and 1 with 1 (we can suppose that the antigen representing string is given in its complementary form). Formally, we define the affinity as the number of identical bits shared by the BCR representing string and $\overline{\mathbf{x}}$.

**Definition 3** For all $\mathbf{x}_i \in \{0,1\}^N$, its affinity with $\overline{\mathbf{x}}$, $aff(\mathbf{x}_i, \overline{\mathbf{x}})$ is given by $aff(\mathbf{x}_i, \overline{\mathbf{x}}) := N - h(\mathbf{x}_i, \overline{\mathbf{x}})$, where $h(\cdot, \cdot) : (\{0,1\}^N \times \{0,1\}^N) \to \{0, \ldots, N\}$ returns the Hamming distance.

**Definition 4** For all $\mathbf{x} = (x_1, \ldots, x_N)$, $\mathbf{y} = (y_1, \ldots, y_N) \in \{0,1\}^N$, their Hamming distance is given by:

$$h(\mathbf{x}, \mathbf{y}) = \sum_{i=1}^{N} \delta_i \qquad \text{where} \qquad \delta_i = \begin{cases} 1 \text{ if } & x_i \neq y_i \\ 0 \text{ otherwise} \end{cases}$$

Other definitions of affinity are often (*e.g.* [37]) constructed as functions of the Hamming distance $aff(\mathbf{x}_i, \overline{\mathbf{x}}) = F(h(\mathbf{x}_i, \overline{\mathbf{x}}))$, for instance with $F$ given by the Gaussian probability density function. These modeling aspect become important when considering the selection mechanism, which is not treated in the present article. Therefore, for our purpose, we can focus on the above definition of affinity.

As first basic mutational rule, we choose to study the one given by single switch-type mutations. At each time step a randomly chosen amino-acid within the BCR binary string switches its amino-acid class. This clearly leads us to a Simple Random Walk (SRW) on $\mathcal{H}_N$. Indeed, we can formalize it as follows:

**Definition 5** Let $\mathbf{X}_n \in \mathcal{H}_N$ be the BCR at step $n$. Let $i \in \{1, \ldots, N\}$ be a randomly chosen index. Then $\mathbf{X}_{n+1} := (X_{n,1}, \ldots, X_{n,i-1}, 1 - X_{n,i}, X_{n,i+1}, \ldots, X_{n,N})$.

*Remark 1* Referring to Definition 2 of neighborhood, as we consider here the standard $N$-dimensional hypercube, $\forall \mathbf{x}_i, \mathbf{x}_j \in \mathcal{H}_N$, $\mathbf{x}_i \sim \mathbf{x}_j \Leftrightarrow h(\mathbf{x}_i, \mathbf{x}_j) = 1$.

We denote the transition probability matrix of the SRW on $\mathcal{H}_N$ by $\mathcal{P}_N$ or simply by $\mathcal{P}$ if no misunderstanding is possible. For all $\mathbf{x}_i, \mathbf{x}_j \in \mathcal{H}_N$:

$$\mathbb{P}(\mathbf{X}_n = \mathbf{x}_j \,|\, \mathbf{X}_{n-1} = \mathbf{x}_i) =: p(\mathbf{x}_i, \mathbf{x}_j) = \begin{cases} 1/N & \text{if } \mathbf{x}_j \sim \mathbf{x}_i \\ 0 & \text{otherwise} \end{cases} \;;\quad \mathcal{P} = (p(\mathbf{x}_i, \mathbf{x}_j))_{\mathbf{x}_i, \mathbf{x}_j \in \mathcal{H}_N}$$

The unique stationary distribution for $\mathcal{P}$ is the homogeneous probability distribution on $\mathcal{H}_N$, denoted by $\boldsymbol{\pi}$: $\forall \mathbf{x}_i \in \mathcal{H}_N$, $\pi_i := \boldsymbol{\pi}(\mathbf{x}_i) = 2^{-N}$. Indeed, $(\mathbf{X}_n)_{n \geq 0}$ is clearly reversible with respect to $\boldsymbol{\pi}$. The uniqueness follows by the Ergodic Theorem.

We also recall a property of $\mathcal{H}_N$ that we will have to deal with: the bipartiteness.

**Definition 6** A graph $G = (V, E)$ is bipartite if there exists a partition of the vertex set $V = V_1 \sqcup V_2$, s.t. every edge connects a vertex in $V_1$ to a vertex in $V_2$.

Typically a bipartition of the hypercube can be obtained by separating the vertices with an odd number of 1 is in their string from those with an even number of 1 is. In Figure 1 we emphasize the bipartite structure of the hypercube $\mathcal{H}_3$.

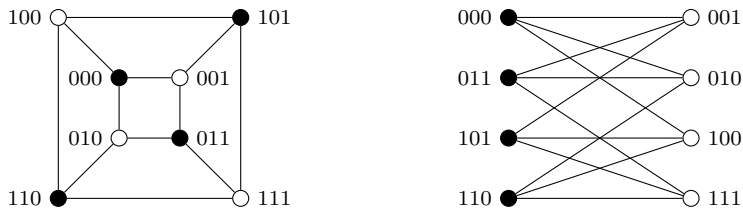

**Figure 1:** Hypercube for $N = 3$ showing its bipartite structure.

A direct and elementary consequence of this property is the periodic behavior of the SRW on $\mathcal{H}_N$, which in particular causes some problems for the convergence through $\boldsymbol{\pi}$. This problem is classically overcome by adding $N$ loops at each vertex, that makes this RW become a *lazy Markov chain* [34]. The corresponding transition probability matrix will be given by $\mathcal{P}_L := (\mathcal{P} + I_{2^N})/2$, where $I_n$ denotes the $n$-dimensional identity matrix.

2.1 Spectral analysis

Most matrices describing the characteristics of the SRW on $\mathcal{H}_N$ can be obtained recursively, thanks to the recursive construction of the hypercube and the operation of cartesian product between two graphs.

**Definition 7** Given two graphs $G_1 = (V_1, E_1)$ and $G_2 = (V_2, E_2)$, the cartesian product between $G_1$ and $G_2$, $G_1 \times G_2$, is a graph with vertex set $V = V_1 \times V_2 = \{(u,v) \,|\, u \in V_1, v \in V_2\}$. Two different vertices $(u_1, v_1)$ and $(u_2, v_2)$ are adjacent in $G_1 \times G_2$ if either $u_1 = u_2$ and $v_1 v_2 \in E_2$ or $v_1 = v_2$ and $u_1 u_2 \in E_1$.

Its a known result [25] that for $N > 1$, $\mathcal{H}_N$ is obtained from $\mathcal{H}_{N-1}$ as: $\mathcal{H}_N = \mathcal{H}_{N-1} \times \mathcal{H}_1$. This characteristic implies the recursive construction of the adjacency matrix and allows to determine the corresponding eigenvalues and eigenvectors. We denote by $A_N$ the adjacency matrix corresponding to $\mathcal{H}_N$; by $I_n$ the $n$-dimensional identity matrix. Then we have:

$$A_1 = \begin{matrix} 0 \\ 1 \end{matrix} \begin{pmatrix} 0 & 1 \\ 1 & 0 \end{pmatrix}; \quad A_2 = \begin{matrix} 00 \\ 01 \\ 10 \\ 11 \end{matrix} \begin{pmatrix} 0 & 1 & 1 & 0 \\ 1 & 0 & 0 & 1 \\ 1 & 0 & 0 & 1 \\ 0 & 1 & 1 & 0 \end{pmatrix} = \begin{pmatrix} A_1 & I_2 \\ I_2 & A_1 \end{pmatrix}$$

Here we wrote in gray the strings corresponding to each row: in order to obtain the adjacency matrices in this form, we simply have to order vertices of $\mathcal{H}_N$ in lexicographical order.

By iteration we obtain [18]:

$$A_n = \begin{pmatrix} A_{n-1} & I_{2^{n-1}} \\ I_{2^{n-1}} & A_{n-1} \end{pmatrix}$$

This iterative construction allows also to determine recursively the spectra of $A_N$ and, consequently, of $\mathcal{P}_N = A_N/N$ (as $\mathcal{H}_N$ is a $N$-regular graph, the transition probability matrix corresponds to the adjacency matrix divided by $N$). Here below we recall the explicit values of the eigenvalues of $A_N$ and $\mathcal{P}_N$ respectively. An extensive proof can be found in [18].

**Theorem 1** *The eigenvalues of $A_N$ are: $N, N-2, N-4, \ldots, -N+4, -N+2, -N$. If we order the $N+1$ distinct eigenvalues of $A_N$ as $\lambda_1^A > \lambda_2^A > \cdots > \lambda_{N+1}^A$, then the multiplicity of $\lambda_k^A$ is $\binom{N}{k-1}$, $1 \leq k \leq N+1$*

**Corollary 1** *The eigenvalues of $\mathcal{P}_N$ are: $1, 1-2/N, 1-4/N, \ldots, -1+4/N, -1+2/N, -1$. If we order the $N+1$ distinct eigenvalues of $\mathcal{P}$ as $\lambda_1 > \lambda_2 > \cdots > \lambda_{N+1}$, then the multiplicity of $\lambda_k$ is $\binom{N}{k-1}$, $1 \leq k \leq N+1$*

Finally we recall the expression of the eigenvectors of $A_N$ (and then also of $\mathcal{P}$), that we will gather together into a matrix. The eigenvectors for $A_1$ are:

$$\mathbf{z}_1 = \begin{bmatrix} 1 \\ 1 \end{bmatrix} \text{ for } \lambda_1^A = 1 \quad \text{and} \quad \mathbf{z}_2 = \begin{bmatrix} 1 \\ -1 \end{bmatrix} \text{ for } \lambda_2^A = -1 \Rightarrow \mathcal{Z}_1 = [\mathbf{z}_1, \mathbf{z}_2]$$

Thanks to the relations between the cartesian product of two graphs and their eigenvectors, it follows by induction that [18]:

$$\mathcal{Z}_n = \left( \begin{array}{c|c} \mathcal{Z}_{n-1} & \mathcal{Z}_{n-1} \\ \hline \mathcal{Z}_{n-1} & -\mathcal{Z}_{n-1} \end{array} \right)$$

Finally, one renormalizes each vector $\mathbf{z}_i$ multiplying it by $\sqrt{2^{-N}}$. We denote by $Q_N$ the resulting matrix, where each column is a $2^N$ vector $\mathbf{v}_i = \sqrt{2^{-N}}\mathbf{z}_i$.

2.2 Evolution of Hamming distances to a fixed node

In this section we focus our attention on the distance process, which is the process obtained from the SRW on $\mathcal{H}_N$ by looking at the Hamming distance between the B-cell representing string at each mutation step and the antigen target representing string. More precisely, $(D_n)_{n \geq 0} := (h(\mathbf{X}_n, \bar{\mathbf{x}}))_{n \geq 0}$ is a RW on $\{0, \ldots, N\}$. From a biological point of view this process represents the evolution of the affinity of our mutating B-cell to the presented antigen. The idea of analyzing the distance of a RW on a graph to some position, where distance means the minimal number of steps that separate two positions, is not unusual. N. Berestycki in [5] applied that to genome rearrangements, where the distance on the graph corresponds biologically to the minimal number of reversals or other mutations needed to transform one genome into the other. Due to the perfect symmetry of the graph we are taking into account and our particular choice of the affinity (which is directly related to the Hamming distance), by studying $(D_n)$ we reduce considerably the number of vertices of the graph, passing from $2^N$ to $N+1$ nodes, without losing the most important properties of the corresponding transition matrix. However, if we consider more complicated models of mutation, it is not possible to reduce the study of the process to the distances to a fixed node. In Figure 2 we show explicitly how pass from $(\mathbf{X}_n)$ to $(D_n)$: since $\bar{\mathbf{x}}$ is fixed and known, we are able to group the vertices by their Hamming distance to $\bar{\mathbf{x}}$. Moreover we keep the original probability of going to the next distance class by considering weighted and directed edges.

The transition probability matrix for $(D_n)$, denoted by $\mathcal{Q}$, is given by Proposition 1 below.

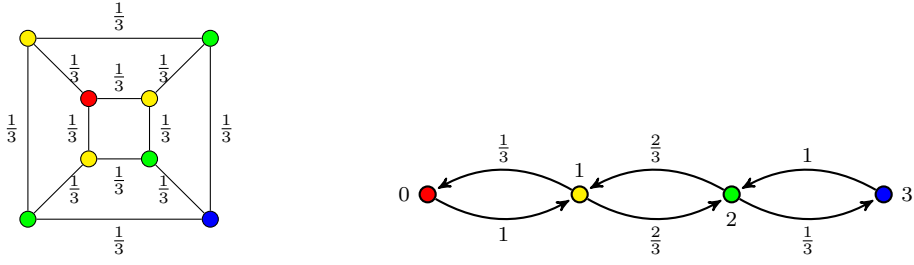

**Figure 2:** From the $(\mathbf{X}_n)$ process (on the left) to the $(D_n)$ process (on the right) (case $N=3$). Near each arrow the probability to travel in the corresponding direction is exhibited. The red vertex always corresponds to $\overline{\mathbf{x}}$, while we represent vertices at the same distance with the same color (yellow for $h=1$, green for $h=2$, and blue for $h=3$).

**Proposition 1** *For all $d, d' \in \{0, \ldots, N\}$:*

$$\mathbb{P}(D_n = d' \mid D_{n-1} = d) =: q(d, d') = \begin{cases} d/N & \text{if } d' = d-1 \\ (N-d)/N & \text{if } d' = d+1 \\ 0 & \text{if } |d'-d| \neq 1 \end{cases} \quad (1)$$

$\mathcal{Q} = (q(d, d'))_{d, d' \in \{0, \ldots, N\}}$ is a $(N+1) \times (N+1)$ tridiagonal matrix where the main diagonal consists of zeros. The stationary distribution for $\mathcal{Q}$ is the binomial probability distribution $\mathcal{B}\left(N, \frac{1}{2}\right) = \left(C_N^d \frac{1}{2^N}\right)_{d \in \{0, \ldots, N\}}$, where $C_N^d = \binom{N}{d} = \frac{N!}{d!(N-d)!}$ is the binomial coefficient. It is the unique stationary distribution for $\mathcal{Q}$: a simple calculation points out the fact that $(D_n)_{n \geq 0}$ is reversible with respect to $\mathcal{B}\left(N, \frac{1}{2}\right)$, then the uniqueness follows again by the Ergodic Theorem.

Anew, we have to deal with bipartiteness: the graph we are taking into account in this section is clearly bipartite, since we can separate its vertices into two subsets containing odd and even nodes respectively and there are no edges connecting two vertices in the same subset. In order to overcome this problem we add $N$ loops at each vertex $\mathbf{x}_i \in \mathcal{H}_N$ which means that the new transition probability matrix for the $(D_n)$ process is, for all $d, d' \in \{0, \ldots, N\}$:

$$\mathbb{P}(D_n = d' \mid D_{n-1} = d) =: q_L(d, d') = \begin{cases} 1/2 & \text{if } d' = d \\ d/(2N) & \text{if } d' = d-1 \\ (N-d)/(2N) & \text{if } d' = d+1 \\ 0 & \text{if } |d'-d| \neq 1 \end{cases} \quad (2)$$

We denote by $\mathcal{Q}_L$ the matrix $\mathcal{Q}_L := (q_L(d, d'))_{d, d' \in \{0, \ldots, N\}}$. Then:

**Proposition 2** $(D_n)_{n\geq 0}$ *converges in law to a binomial random variable with parameters $N$ and $1/2$. Explicitly:*

$$(\mathcal{Q}_L)_d \to \mathcal{B}\left(N, \frac{1}{2}\right)_d \quad for \quad n \to +\infty$$

*Proof* The proof follows directly observing that $\mathcal{Q}_L$ represents an irreducible and, now, aperiodic MC, with the same stationary distribution as $\mathcal{Q}$ (see [41] for a proof of the general result). □

The spectral analysis of $\mathcal{Q}$ gives us the following result.

**Theorem 2** *For fixed $N$, the spectra of the transition probability matrix $\mathcal{Q}$ corresponding to the $(D_n)$ process is composed by the same $N+1$ distinct eigenvalues as the spectra of $\mathcal{P}$, each with multiplicity 1.*

*Proof* The proof consists of a simple calculation of the eigenvalues of matrix $\mathcal{Q}$ for little $N$. Then we reason by iteration. We can also give the system we use for determining the eigenvectors. For fixed $N$ let us denote by $\lambda_{\pm k}$ the eigenvalue $\frac{\pm(N-2k)}{N}$ for $0 \leq k \leq \lfloor N/2 \rfloor$. We denote by $\mathbf{x}_{\pm k}$ the corresponding unknown eigenvector. Then we have the following matrix equation:

$$\mathcal{Q}\mathbf{x}_{\pm k} = \lambda_{\pm k}\mathbf{x}_{\pm k}$$

Which is:

$$\begin{cases} x_{\pm k,2} = \lambda_{\pm k} x_{\pm k,1} \\ \\ \frac{1}{N}x_{\pm k,1} + \frac{N-1}{N}x_{\pm k,3} = \lambda_{\pm k} x_{\pm k,2} \\ \\ \frac{2}{N}x_{\pm k,2} + \frac{N-2}{N}x_{\pm k,4} = \lambda_{\pm k} x_{\pm k,3} \\ \\ \quad\quad\quad\quad\quad \vdots \\ \\ \frac{N-1}{N}x_{\pm k,N-1} + \frac{1}{N}x_{\pm k,N+1} = \lambda_{\pm k} x_{\pm k,N} \\ \\ x_{\pm k,N} = \lambda_{\pm k} x_{\pm k,N+1} \end{cases}$$

□

*Remark 2* Using the classical results of S. N. Ethier and T. G. Kurtz [14] it is possible to prove that, denoting by $x_N(t)$ the process $x_N(t) = \frac{D_{\lfloor Nt \rfloor}}{N}$, it converges in probability through $x(t)$, solution of the differential equation $\dot{x}(t) = -2x(t) + 1$ on a finite time window:

$$\forall \varepsilon > 0, \forall T > 0, \mathbb{P}\left(\sup_{t \in [0,T]} |x_N(t) - x(t)| > \varepsilon\right) \to 0 \quad \text{for } N \to \infty.$$

*Remark 3* We can easily observe that $x(t)$ rapidly converges to $1/2$ for all $x_0 \in [0,1]$. In particular if we start at $x_0 = 1/2$, we stay there for all $t$. That suggests that the $(D_n)$ process, for $N$ going to infinity, reaches a value of about $N/2$ exponentially fast, and then tends to remain there.

From an heuristic viewpoint we can explain how we derived the above equation. First of all, we take into account the following rescaled process:

$$x_n := D_n/N$$

As $(D_n) \in \{0, \ldots, N\}$, $x_n \in [0,1]$. Denoting by $q_n(x) = \mathbb{P}(x_n = x)$ and using Equation (1), we have:

$$q_{n+1}(x) = (1-x)q_n\left(x - \frac{1}{N}\right) + xq_n\left(x + \frac{1}{N}\right)$$

Now we apply the Taylor is theorem for $N \gg 1$:

$$q_{n+1}(x) = (1-x)\left(q_n(x) - \frac{1}{N}q'_n(x) + o\left(\frac{1}{N}\right)\right) + x\left(q_n(x) + \frac{1}{N}q'_n(x) + o\left(\frac{1}{N}\right)\right)$$

From which we get:

$$q_{n+1}(x) - q_n(x) = \frac{1}{N}(x - (1-x))q'_n(x) + o\left(\frac{1}{N}\right)$$

Defining the process $\tilde{q}(t,x) = q_{\lfloor Nt \rfloor}(x)$, with $t = \frac{n}{N}$, we obtain:

$$\partial_t \tilde{q}(t,x) = (2x-1)\partial_x \tilde{q}(t,x) + o\left(\frac{1}{N}\right)$$

And consequently, the corresponding transport equation is:

$$\partial_t q(t,x) = (2x-1)\partial_x q(t,x) \qquad (3)$$

The differential equation associated with Equation (3) (its characteristic equation) is:

$$\dot{x}(t) = -2x(t) + 1$$

which has solution:

$$x(t) = \frac{1}{2} + \left(x_0 - \frac{1}{2}\right)e^{-2t}$$

It is also possible to derive a diffusion approximation by expanding the generator at second order.

2.3 Hitting times

In this section we give explicit formulas to compute the hitting time from node $\mathbf{x}_i$ to $\mathbf{x}_j$: the expected number of steps before $\mathbf{x}_j$ is visited, starting from $\mathbf{x}_i$. More precisely, we define by $\tau_{\{\mathbf{x}_j\}} := \inf\{n \geq 0 \,|\, \mathbf{X}_n = \mathbf{x}_j\}$: we are interested in studying its expectation, $\mathbb{E}_{\mathbf{x}_i}[\tau_{\{\mathbf{x}_j\}}]$. The formula we gave in Section 2.3.1 is directly obtained from the more general one given by L. Lovász in [36] (we recall Equation (6) simply because we will need it later). On the other hand, the formula given in Section 2.3.2 is obtained from the $(D_n)$ process and the procedure is inspired by those one used in [33].

### 2.3.1 Analysis of $\mathbb{E}_{\mathbf{x}_0}[\tau_{\{\overline{\mathbf{x}}\}}]$ using the spectrum of $\mathcal{P}$.

**Definition 8** Let $H$ be the $2^N \times 2^N$ symmetric matrix having as $(i,j)^{\text{th}}$ entry: $(H)_{ij} = H(i,j) = \mathbb{E}_{\mathbf{x}_i}[\tau_{\{\mathbf{x}_j\}}]$ for all $\mathbf{x}_i, \mathbf{x}_j \in \mathcal{H}_N$. Clearly $H(i,i) = 0$ for all $i$.

The $N$-regularity of the graph implies that:

$$H(i,j) = 1 + \sum_{\{k\,|\,h(i,k)=1\}} \mathcal{P}_{ik} H(k,j) = 1 + \frac{1}{N} \sum_{\{k\,|\,h(i,k)=1\}} H(k,j) \quad \text{for } i \neq j \quad (4)$$

To relate the hitting time with the spectrum, we first define $F := J + \mathcal{P}H - H$, where $J$ is a $2^N \times 2^N$ matrix whose entries are all 1. From Equation (4), it follows that $F$ is a diagonal matrix, as $(H)_{ij} = (J)_{ij} + (\mathcal{P}H)_{ij}$ for $i \neq j$. Moreover $F'\boldsymbol{\pi} = \mathbf{1}$, where $\mathbf{1} = (1,\ldots,1)'$, since

$$F'\boldsymbol{\pi} = (J + (\mathcal{P} - I_{2^N})H)'\boldsymbol{\pi} = J\boldsymbol{\pi} + H'(\mathcal{P} - I_{2^N})'\boldsymbol{\pi} = J\boldsymbol{\pi} + H'(\mathcal{P}'\boldsymbol{\pi} - \boldsymbol{\pi}) = J\boldsymbol{\pi} = \mathbf{1}$$

Therefore, we deduce that $F = 2^N I_{2^N}$ and that $H$ is solution of

$$(I_{2^N} - \mathcal{P})H = J - 2^N I_{2^N} \tag{5}$$

As shown in [36], the solution is given by $H := (I_{2^N} - \mathcal{P} + \mathbf{1}\boldsymbol{\pi}')^{-1}(J - 2^N I_{2^N})$, yielding:

**Theorem 3** *Given a SRW on $\mathcal{H}_N$, the hitting time from vertex $i$ to $j$ is given by:*

$$H(i,j) = 2^N \sum_{k=2}^{2^N} \frac{1}{1-\lambda_k}(v_{kj}^2 - v_{ki}v_{kj}), \tag{6}$$

*where $\lambda_k$ is the $k^{\text{th}}$ eigenvalue of $\mathcal{P}$ and $v_{ki}$ corresponds to the $i^{\text{th}}$ component of the $k^{\text{th}}$ eigenvector of $\mathcal{P}$, as given in Section 2.1.*

*2.3.2 Analysis of $\mathbb{E}_{\mathbf{x}_0}[\tau_{\{\overline{\mathbf{x}}\}}]$ from the $D_n$ viewpoint.*

For the sake of simplicity, we denote $H(D_0) := \mathbb{E}_{\mathbf{x}_0}[\tau_{\{\overline{\mathbf{x}}\}}]$ as it depends only on the initial Hamming distance of $\mathbf{X}_0$ to $\overline{\mathbf{x}}$, $D_0$.

*Remark 4* Due to (1), starting at point $\mathbf{x}_0$ with $D_0 = \overline{d}$, we have:

$$\begin{cases} \mathbb{P}(D_1 = \overline{d}+1 \,|\, D_0 = \overline{d}) =: q(\overline{d}, \overline{d}+1) = (N - \overline{d})/N \\ \mathbb{P}(D_1 = \overline{d}-1 \,|\, D_0 = \overline{d}) =: q(\overline{d}, \overline{d}-1) = \overline{d}/N \end{cases}$$

We are now able to define a new recursive formula for (4), which will be more convenient if evaluated explicitly:

$$H(\overline{d}) = 1 + \frac{N - \overline{d}}{N} H(\overline{d}+1) + \frac{\overline{d}}{N} H(\overline{d}-1) \tag{7}$$

with boundary conditions:

$$H(0) = 0 \text{ and } H(1) = 2^N - 1 = \sum_{j=0}^{N} C_N^j - 1 \tag{8}$$

Taking the difference $\Delta(\overline{d}) := H(\overline{d}) - H(\overline{d}-1)$, we obtain:

$$\Delta(\overline{d}+1) = H(\overline{d}+1) - H(\overline{d}) = \frac{\overline{d}}{N}\left(\Delta(\overline{d}+1) + \Delta(\overline{d})\right) - 1$$

And finally:

$$\Delta(\overline{d}+1) = \frac{\overline{d}}{N - \overline{d}} \Delta(\overline{d}) - \frac{N}{N - \overline{d}} \quad \text{with } \Delta(1) = H(1) \tag{9}$$

Then we can prove by iteration the following result:

**Theorem 4** *Given a SRW on $\mathcal{H}_N$, the hitting time to cover a Hamming distance equal to $\overline{d}$, $H(\overline{d})$ with $0 \leq \overline{d} \leq N$ is obtained as:*

$$H(\overline{d}) = \sum_{d=0}^{\overline{d}-1} \frac{\sum_{j=1}^{N-1-d} C_N^{d+j} + 1}{C_{N-1}^d} \tag{10}$$

*Proof* One have to prove that:

$$\Delta(d+1) = \frac{\sum_{j=1}^{N-1-d} C_N^{d+j} + 1}{C_{N-1}^d} \tag{11}$$

$$\Delta(d+1) = \frac{d \cdot \Delta(d)}{N-d} - \frac{N}{N-d} = \frac{d}{N-d}\left(\frac{(d-1) \cdot \Delta(d-1)}{N-(d-1)} - \frac{N}{N-(d-1)}\right) - \frac{N}{N-d}$$

$$= \frac{d(d-1) \cdot \Delta(d-1)}{(N-d)(N-(d-1))} - N\left(\frac{d}{(N-d)(N-(d-1))} + \frac{1}{N-d}\right) \tag{12}$$

Proceeding by iteration we obtain two terms, where the first one multiplies $\Delta(1)$. From Equation (9) we know that $\Delta(1) = H(1) = \sum_{j=0}^{N} C_N^j - 1$. A convenient use of the properties of the factorial operator allows us to reach the following expression:

$$
\begin{aligned}
(12) &= \frac{d!(N-1-d)!}{(N-1)!}\left(\sum_{j=0}^{N} C_N^j - 1\right) - N\left(\frac{d!(N-1-d)!}{(N-1)!} + \frac{d!(N-1-d)!}{2!(N-2)!} + \cdots \right.\\
&\quad \left. + \frac{d!(N-1-d)!}{(d-1)!(N-(d-1))!} + \frac{d!(N-1-d)!}{d!(N-d)!}\right) = \\
&= \frac{d!(N-1-d)!}{(N-1)!}\left(1 + \sum_{j=1}^{N-1-d} \frac{N!}{(d+j)!(N-(d+j))!}\right) = \frac{\sum_{j=1}^{N-1-d} C_N^{d+j} + 1}{C_{N-1}^d}
\end{aligned}
$$

By using again (9), we can now easily express $H(\overline{d})$ in the following way

$$H(\overline{d}) = \sum_{d=0}^{\overline{d}-1} \Delta(d+1) = \sum_{d=0}^{\overline{d}-1} \frac{\sum_{j=1}^{N-1-d} C_N^{d+j} + 1}{C_{N-1}^d}$$

which can be evaluated for reasonable values of $N$.  □

We can immediately observe that $H(\overline{d})$ is a monotonically increasing function. Moreover, $H$ is concave. Indeed, thanks to Proposition 4 we can prove that $\forall d \in \{1, \ldots, N-1\}$:

$$H(d) - H(d-1) \geq H(d+1) - H(d) \iff \Delta(d) \geq \Delta(d+1)$$

Furthermore, we can evaluate the following limit:

$$\lim_{N \to \infty} \frac{H(\alpha N)}{2^N} \quad \text{for } \alpha \in \,]0, 1]. \tag{13}$$

*Remark 5* The case $\alpha = 0$ is trivial: if $\alpha = 0$ this limit is equal to 0 since $H(0) = 0$.

*Remark 6* Proposition 3 below, which evaluates (13), confirms the statement made in Remark 3: as $N$ goes to infinity, $(D_n)$ goes quickly to $N/2$ and then $H(d)$ is always of order $\sim 2^N$ irrespective of $d \neq 0$.

**Proposition 3** *For all $\alpha \in \,]0, 1]$:*

$$\lim_{N \to \infty} \frac{H(\alpha N)}{2^N} = 1$$

*Proof* Since $H$ is an increasing function and by using Equation (10) we have:

$$2^N - 1 = H(1) \leq H(\alpha N) \leq H(N) = \sum_{d=0}^{N-1}\sum_{j=1}^{N-1-d} \frac{C_N^{d+j}}{C_{N-1}^d} + \sum_{d=0}^{N-1} \frac{1}{C_{N-1}^d} =: I_1 + I_2$$

We examine the two terms of the last member separately.

$$I_2 \leq 2 + \frac{2}{N-1} + (N-4)\frac{2}{(N-1)(N-2)} \tag{14}$$

We can prove it just by looking at Pascal is triangle.

Now, if we consider $I_1$, we see that there is no contribution for $d = N-1$, as the internal sum is zero valued. Moreover we have:

$$\sum_{j=1}^{N-1-d} C_N^{d+j} \leq \sum_{j=0}^{N} C_N^j = 2^N$$

And so:

$$I_1 \leq 2^N \sum_{d=0}^{N-2} \frac{1}{C_{N-1}^d} \stackrel{(14)}{\leq} 2^N \left(1 + \frac{2}{N-1} + (N-4)\frac{2}{(N-1)(N-2)}\right)$$

By putting together all these inequalities and dividing by factor $2^N$ we get that:

$$1 - \frac{1}{2^N} \leq \frac{H(\alpha N)}{2^N} \leq 1 + \frac{2}{N-1} + \frac{2(N-4)}{(N-1)(N-2)} + \frac{1}{2^N}\left(2 + \frac{2}{N-1} + \frac{2(N-4)}{(N-1)(N-2)}\right)$$

The result comes directly by applying the squeeze theorem. □

This result can be extended to a SRW on a generic state-space $\mathcal{S}^N$, with $|\mathcal{S}| = s$. More precisely, one can prove in a similar way as we did for $\mathcal{H}_N$ the following result:

**Proposition 4** *The order of magnitude of the hitting time for a switch-type mutational model on the state-space $\mathcal{S}^N$, with $|\mathcal{S}| = s$, is $s^N$, for $N$ big enough.*

This is the consequence of Theorem 5 and Proposition 5 below.

**Theorem 5** *Given a SRW on $S^N$, the hitting time to cover a Hamming distance equal to $\bar{d}$, $H_s(\bar{d})$ with $0 \leq \bar{d} \leq N$ is obtained as:*

$$H_s(\bar{d}) = \sum_{d=0}^{\bar{d}-1} \frac{\sum_{j=d+1}^{N} C_N^j (s-1)^j}{C_{N-1}^d (s-1)^d} \tag{15}$$

**Proposition 5** *For all $\alpha \in ]0,1]$:*

$$\lim_{N \to \infty} \frac{H_s(\alpha N)}{s^N} = 1$$

## 3 More mutational models: how does the structure of the hypercube change?

In this section, we explore other mutation rules, which change the internal graph structure of the hypercube, therefore the dynamics of the RW and the characteristic time-scales of the exploration of the state-space.

3.1 Study of various mutation rules

In this section, we study mainly three mutation rules : a model of switch of $k$-length strings, a model of permutation of two bits and a model of switch of 1 or 2-length strings depending on the Hamming distance to a fixed node representing the antigen target cell.

*3.1.1 The exchange mutation model.*

We consider a model where given an initial B-cell representing string, each mutation step consists in permuting two randomly chosen bits.

**Definition 9** Let $\mathbf{X}_n \in \{0,1\}^N$ be the BCR at step $n$. Let $i \in \{1,\ldots,N\}$, $j \in \{1,\ldots,N\} \setminus \{i\}$ two randomly chosen indexes. Then (we can suppose, without loss of generality, that $j > i$):

$$\mathbf{X}_{n+1} = (X_{n,1},\ldots,X_{n,i-1},X_{n,j},X_{n,i+1},\ldots,X_{n,j-1},X_{n,i},X_{n,j+1},\ldots,X_{n,N})$$

With this mutation rule, we loose a very important property : the connectivity of the graph. We denote by $\mathcal{H}_{(s)} \subset \{0,1\}^N$ the set containing the $C_N^s$ vertices having $s$ 1 in their strings. The state-space $\{0,1\}^N$ is divided into $N+1$ connected components: $\mathcal{H}_{(s)}$, $0 \leq s \leq N$.

**Proposition 6** *There are exactly $\frac{N(N-1)}{2}$ (non-oriented) edges ending at each vertex counting the possible loops. Each node $\mathbf{x} \in \mathcal{H}_{(s)}$ has exactly $\frac{(N-s)^2-(N-s^2)}{2}$ loops.*

**Corollary 2** $\mathbb{P}(\mathbf{X}_n = \mathbf{x}_j | \mathbf{X}_{n-1} = \mathbf{x}_j) = \frac{(N-s)^2-(N-s^2)}{N(N-1)}$. *Then the probability of remaining on the same node is 1 if $s=0$ or $s=N$.*

*Proof (Proposition 6)* The first statement is obtained by simple combinatory arguments. Let us consider $\mathbf{x} \in \mathcal{H}_{(s)}$ with $0 \leq s \leq N$: it is composed by exactly $s$ ones and $N-s$ zeros. For the sake of clarity let us consider that $\{0,\ldots,N\} = I \sqcup J$ so that $|I| = s$, $|J| = N-s$ and $x_i = 1 \; \forall i \in I$, $x_j = 0 \; \forall j \in J$. We obtain a loop each time we choose both random indices either in $I$ ($C_s^2$ possibilities) or in $J$ ($C_{N-s}^2$ possibilities). Then the total number of loops is obtained by the sum of these two cases, *i.e.* $\frac{(N-s)^2-(N-s^2)}{2}$. □

We can also describe qualitatively the behavior of the $(D_n)$ process referring to this current model. As a general principle, we have that $D_n = D_{n-1} + i$, $i \in \{0, \pm 2\}$. Therefore, clearly $\mathbb{P}(D_n = d'|D_{n-1} = d) = 0$ if $|d' - d| > 2$ or $|d' - d| = 1$. Moreover, we have maximal and minimal values of $D_n$ depending on $s_0$ and $\overline{s}$ so that $\mathbf{X}_0 \in \mathcal{H}_{(s_0)}$ and $\overline{\mathbf{x}} \in \mathcal{H}_{(\overline{s})}$. Indeed:

**Proposition 7** *Given $\overline{\mathbf{x}} \in \mathcal{H}_{(\overline{s})}$ and $\mathbf{X}_0 \in \mathcal{H}_{(s_0)}$, then $\forall n \geq 0$:*

$$\begin{cases} |\overline{s} - s_0| \leq D_n \leq \overline{s} + s_0 & \text{if } \overline{s} + s_0 \leq N \\ \\ |\overline{s} - s_0| \leq D_n \leq (N - \overline{s}) + (N - s_0) & \text{if } \overline{s} + s_0 > N \end{cases}$$

*Proof* The proof follows immediately by counting how many possibilities there are to arrange $s$ ones and $N - s$ zeros in a $N$-length string. □

*Remark 7* From Proposition 7 one can see that if $\overline{s} = s_0 =: s$ and $2s \neq N$ then:

$$0 \leq D_n < N$$

*3.1.2 Class switch of k-length strings.*

Let $\mathbf{X}_0 = (X_{0,1}, \ldots, X_{0,N}) \in \{0,1\}^N$ be the B-cell entering the somatic hypermutation process. At each mutation step we switch the class of $k$ consecutive amino-acids.

**Definition 10** Let $\mathbf{X}_n \in \{0,1\}^N$ be the BCR at step $n$. Let $i \in \{1, \ldots, N - (k-1)\}$ be a randomly chosen index. Then $\mathbf{X}_{n+1} := (X_{n,1}, \ldots, X_{n,i-1}, 1 - X_{n,i}, \ldots, 1 - X_{n,i+k-1}, X_{n,i+k}, \ldots, X_{n,N})$.

*Remark 8* If $k = 1$ we are in the case of a SRW on $\mathcal{H}_N$.
If $k = N$ we stay on a 2-length cycle. Indeed we have that $\mathbf{X}_l = \mathbf{X}_0$ for $l$ even and $X_l = \mathbf{1} - \mathbf{X}_0$ for $l$ odd. For this reason the case $k = N$ does not appear interesting neither from a mathematical nor from a biological point of view.

Here below we give some basic properties of this RW, that one can easily prove by simple combinatory arguments.

**Proposition 8** *Each vertex has exactly $N - (k-1)$ neighbors and no loops. Therefore, for all $\mathbf{x}_i, \mathbf{x}_j$ in $\{0,1\}^N$:*

$$\mathbb{P}(\mathbf{X}_n = \mathbf{x}_j | \mathbf{X}_{n-1} = \mathbf{x}_i) =: p_k(i,j) = \begin{cases} \dfrac{1}{N - (k-1)} & \text{if } \mathbf{x}_j \sim \mathbf{x}_i \\ \\ 0 & \text{otherwise} \end{cases}$$

*Remark 9* As regards to this current model, given $\mathbf{x}_i, \mathbf{x}_j \in \{0,1\}^N$, we have: $\mathbf{x}_i \sim \mathbf{x}_j \Leftrightarrow h(\mathbf{x}_i, \mathbf{x}_j) = k$ and the $k$ different elements have consecutive indexes.

Thus, $\mathcal{P}_k = (p_k(\mathbf{x}_i, \mathbf{x}_j))_{\mathbf{x}_i, \mathbf{x}_j \in \mathcal{H}_k}$ is the $2^N \times 2^N$ transition probability matrix.

For fixed $k \in \{1, \ldots, N\}$ the graph underlying the RW corresponding to the model of class switch of $k$-length strings has exactly $2^{k-1}$ connected components, each one composed of $2^{N-(k-1)}$ elements.
Because of the non connectivity of the graph, we can focus on the connected component to which $\mathbf{X}_0$ belongs and find out the properties of our RW on it. For fixed $N$ and $k$ and dealing with each connected component separately, we are describing a SRW on a $(N-(k-1))$-hypercube. Henceforth we obtain $2^{k-1}$ distinct hypercube-type structures of the same size.

We can limit our study to the connected component containing $\mathbf{X}_0$, which is, up to a change of variables, a $(N-(k-1))$-dimensional hypercube. Let $\overline{\mathcal{P}}_k$ be the restriction of $\mathcal{P}_k$ to this connected component. If we conveniently order the $2^{N-(k-1)}$ distinct vertices, than $\overline{\mathcal{P}}_k = \mathcal{P}_{N-(k-1)}$. At this stage, it is possible to translate all classical results we know about the SRW on $\mathcal{H}_n$, for $n = N-(k-1)$, on each connected component of this current graph, remembering the definition of neighborhood given in Remark 9.

*3.1.3 Class switch of 1 or 2-length strings depending on the Hamming distance to $\overline{\mathbf{x}}$.*

The models we described in Sections 3.1.1 and 3.1.2 present an important limitation: the underlying graphs are non-connected. Due to the choice we made of affinity, a model which does not enable to explore the whole state-space is not very relevant. Indeed, if the graph is non-connected and the target chain does not belong to the connected component containing the B-cell which first enters the somatic hypermutation process, then we never reach the target configuration. From a biological viewpoint, it may be more relevant to consider a smoother affinity model, in which the BCR representing string reaches the target when most, but not all, bits are similar. In this case, considering a non-connected graph, is not necessarily a problem.

Another way to overcome the problem of non-connectivity is to consider a model which allows to vary the length of the strings submitted to switch-type mutations. Moreover, it is biologically credible that during the GC process B-cells can modify their mutation rate, making it somehow inversely proportional to their affinity to the antigen [45,7,4]: the greater the affinity, the lower is the mutation rate. Indeed, B-cells during the GC process compete for different rescue signals (from Helper T-cells or FDCs), and that determines their fate: undergo further mutations or differentiate into plasma cells or memory cells ([1], Chapters 7). We found the hypothesis that the regulation of the hypermutation process is dependent on receptor affinity also in other works, as [9] by L.N. De Castro and F. J. Von Zuben, where they proposed a computational implementation of the clonal selection principle to design genetic optimization algorithms, taking into account AAM during an adaptive immune response. In terms of our mathematical model, we can translate it by making the size $k$ of the strings which can mutate to be directly proportional to the Hamming distance to $\overline{\mathbf{x}}$ at each mutation step:

$k_n = f(D_n)$, with $f : \{0,\dots,N\} \to \{0,\dots,N\}$ being an increasing function.

Despite many choices of the function $f$ are possible, hereinafter we consider a very elementary case, where $f$ is a step function on two intervals.

**Definition 11** Let $\mathbf{X}_n \in \{0,1\}^N$ be the BCR at step $n$. We denote by $k_n$:

$$k_n := f(D_n) = \begin{cases} 1 \text{ if } D_n \leq 1 \\ 2 \text{ if } D_n > 1 \end{cases}$$

Let $i \in \{1,\dots,N-(k_n-1)\}$ be a randomly chosen index. Then:
$\mathbf{X}_{n+1} := (X_{n,1},\dots,X_{n,i-1}, 1-X_{n,i},\dots, 1-X_{n,i+k_n-1}, X_{n,i+k_n},\dots,X_{n,N})$.

This model is an interesting and simple way to generalize the basic mutational model without loosing the property of connectivity of the graph. The addition of this flexibility was not only motivated by biological reasons, but we also expect that this modification decreases the hitting time to a fixed node. This is actually true: the hitting time is halved compared to the basic model (at least for $N$ big enough). We will also show that the stationary distribution is concentrated on a half part of the hypercube, the one to whom $\overline{\mathbf{x}}$ belongs.

*Remark 10* For fixed $N$ and $k = 2$ the graph is divided into two connected components composed of $2^{N-1}$ vertices. Two nodes belonging to the same connected component have a Hamming distance of $2t$ with $0 \leq t \leq \lfloor N/2 \rfloor$. On the other hand, two vertices belonging to different connected components have a Hamming distance of $(2t+1)$ with $0 \leq t \leq \lfloor (N-1)/2 \rfloor$.

In order to analyze this process, we have to distinguish two cases. For fixed $N$ and $\overline{\mathbf{x}}$, the process we obtain:

**case 1: $D_0 = 2t$, $t > 0$.** $\mathbf{X}_0$ belongs to the same connected component as $\overline{\mathbf{x}}$, so we are working on a $(N-1)$-dimensional hypercube, following the model of class switch of 2-length strings. we stay in this connected component all over the process till we arrive at $\overline{\mathbf{x}}$, as it is impossible to obtain a Hamming distance equal to 1.

**case 2: $D_0 = 2t+1$, $t > 0$.** We necessarily take $k = 2$ and Remark 10 implies that $\mathbf{X}_0$ belongs to a different connected component than $\overline{\mathbf{x}}$. In order to reach the connected component containing $\overline{\mathbf{x}}$, we have to visit a node $\mathbf{x}^*$ so that $h(\mathbf{x}^*, \overline{\mathbf{x}}) = 1$, and $|\{\mathbf{x}^* \,|\, h(\mathbf{x}^*, \overline{\mathbf{x}}) = 1\}| = N$. Then, if $D_0 = 1$ we are allowed to change only one element of the B-cell representing string. With probability $1/N$ we arrive directly at $\overline{\mathbf{x}}$ and with probability $(N-1)/N$ we obtain $D_1 = 2$. Then we go back to case 1.

**Proposition 9** *The graph corresponding to the current model is divided into two connected components: $\mathcal{H}_N^{(1-2)}$ and its complementary $\overline{\mathcal{H}_N}^{(1-2)}$, s.t. $\overline{\mathbf{x}} \in \overline{\mathcal{H}_N}^{(1-2)}$. $\overline{\mathcal{H}_N}^{(1-2)}$ is accessible from $\mathcal{H}_N^{(1-2)}$, but not conversely. Vertices belonging to $\overline{\mathcal{H}_N}^{(1-2)}$ are positive recurrent and vertices belonging to $\mathcal{H}_N^{(1-2)}$ are transient.*

*Proof* The existence of two connected components depends on the use of the model of switch of 2-length strings. Indeed the structure of the graph we are considering here essentially corresponds to that of the graph underlying the model of switch of 2-length strings, up to the addition of some oriented edges from $\mathcal{H}_N^{(1-2)}$ to $\overline{\mathcal{H}_N}^{(1-2)}$. As long as we stay in $\overline{\mathcal{H}_N}^{(1-2)}$ or $\mathcal{H}_N^{(1-2)}$ we are just allowed to switch 2-length strings. Moreover, we have already observed that when we are in $\overline{\mathcal{H}_N}^{(1-2)}$ we can't exit, while when we are in $\mathcal{H}_N^{(1-2)}$ we can reach $\overline{\mathcal{H}_N}^{(1-2)}$ by visiting one among the $N$ nodes having Hamming distance 1 from $\overline{\mathbf{x}}$, and that happens in a finite number of steps. Therefore:

$$\begin{cases} \mathbb{P}(\tau_{\mathbf{x}_i} < \infty) = 1 \text{ for all } \mathbf{x}_i \in \overline{\mathcal{H}_N}^{(1-2)} \Rightarrow \mathbf{x}_i \text{ is recurrent} \\ \\ \mathbb{P}(\tau_{\mathbf{x}_i} < \infty) < 1 \text{ for all } \mathbf{x}_i \in \mathcal{H}_N^{(1-2)} \quad \Rightarrow \mathbf{x}_i \text{ is transient} \end{cases}$$

In particular, vertices belonging to $\overline{\mathcal{H}_N}^{(1-2)}$ are positive recurrent as the chain is irreducible on $\overline{\mathcal{H}_N}^{(1-2)}$ and $|\overline{\mathcal{H}_N}^{(1-2)}| < \infty$. □

The following known result about stochastic processes, justify Corollary 3 below.

**Theorem 6** *Let $(\mathbf{X}_n)_{n \geq 0}$ be a Markov chain on a state-space $\mathcal{S}$ and $\mathbf{x}_i \in \mathcal{S}$ be positive recurrent. Let $m_i$ be the mean return time: $m_i = \mathbb{E}(\tau_{\{\mathbf{x}_i\}} | \mathbf{X}_0 = \mathbf{x}_i)$. Denoting by $\mathcal{S}_r \subseteq \mathcal{S}$ the positive recurrent connected component to which $\mathbf{x}_i$ belongs, then a stationary distribution $\overline{\boldsymbol{\pi}}$ is given by:*

$$\overline{\pi}_i = m_i \quad \forall \mathbf{x}_i \in \mathcal{S}_r$$
$$\overline{\pi}_i = 0 \quad \forall \mathbf{x}_i \in \mathcal{S} \setminus \mathcal{S}_r$$

Theorem 6 is proven by considering the relations among recurrent and transient classes, stationary distributions and return time (see [41] for some more details).

**Corollary 3** *The stationary distribution for the RW we describe in the present section, $\overline{\boldsymbol{\pi}}$, is given by:*

$$\overline{\pi}_i = \begin{cases} \dfrac{1}{2^{N-1}} & \text{if } \mathbf{x}_i \in \overline{\mathcal{H}_N}^{(1-2)} \\ \\ 0 & \text{if } \mathbf{x}_i \in \mathcal{H}_N^{(1-2)} \end{cases} \quad (16)$$

Corollary 3 is a consequence of Theorem 6 and the study of the SRW on an $N$-dimensional hypercube.

*3.1.4 Allowing 1 to k mutations*

In this section we analyze how the *N*-dimensional hypercube changes if we allow 1 to *k* independent switch-type mutations at each step, with *k* fixed, $k \leq N$.

**Definition 12** Let $\mathbf{X}_n \in \{0,1\}^N$ be the BCR at step *n*. Let *k* be an integer, $1 \leq k \leq N$ and $\forall i$, $1 \leq i \leq k$, $a_i := \mathbb{P}(i \text{ independent switch-type mutations})$. Then with probability $a_i$, $\mathbf{X}_{n+1}$ is obtained from $\mathbf{X}_n$ by repeating *i* times, independently, the process described by Definition 5.

By definition, the corresponding transition probability matrix is a convex combination of $\mathcal{P}^i$, for $1 \leq i \leq k$ ($\mathcal{P}^i$ is the transition probability matrix corresponding to *i* iterations of the process of a single bit mutation):

$$\sum_{i=1}^{k} a_i \mathcal{P}^i, \quad \text{with } \sum_{i=1}^{k} a_i = 1. \tag{17}$$

**Definition 13** Let us fix $a_i = 1/k$ $\forall i$. We denote by $\mathcal{P}^{(k)} := 1/k \sum_{i=1}^{k} \mathcal{P}^i$. Accordingly, we denote the graph underlying this RW $\mathcal{H}_N^{(k)}$.

*Remark 11* Since the mutations are assumed to be independent, then *k* represents the maximum Hamming distance the process can cover in a single mutation step. Thanks to the independence of each single mutation, two or more mutations may nullify their respective action: in particular for $k \geq 2$ there is a non-zero probability of remaining at the same position. From a biological point of view, this behavior can be interpreted as the possibility of doing mutations which have no effect on the BCR structure.

We can now evaluate the eigenvalues of $\mathcal{P}^{(k)}$, $\lambda_j^{(k)}$ by using the eigenvalues $\lambda_j$ of $\mathcal{P}$ (Section 2.1). Due to the fact that all $\mathcal{P}^i$ commute with each other, the eigenvalues are given by:

$$\lambda_j^{(k)} = \frac{1}{k} \sum_{i=1}^{k} \lambda_j^i \tag{18}$$

and $\mathcal{P}^{(k)}$ and $\mathcal{P}$ have the same eigenvectors. We give explicitly the expression of all $\lambda_i^{(k)}$ and concentrate on the second largest eigenvalue, $\lambda_2^{(k)}$.

**Proposition 10** *The $N+1$ distinct eigenvalues of matrix $\mathcal{P}^{(k)}$ are:*

- $\lambda_1^{(k)} = 1$ ;
- $\lambda_j^{(k)} = \dfrac{\lambda_j}{k} \cdot \dfrac{1 - \lambda_j^k}{1 - \lambda_j}$ *for $2 \leq j \leq N$ ;*
- $\lambda_{N+1}^{(k)} = \dfrac{1}{2k}\left((-1)^k - 1\right) = \begin{cases} 0 & \text{if } k \text{ is even} \\ -1/k & \text{if } k \text{ is odd} \end{cases}$

*The multiplicity of $\lambda_j^{(k)}$ is $\binom{N}{j-1}$, $1 \leq j \leq N+1$*

*Proof* This result comes directly from the evaluation of Equation (18), for the already known values of all $\lambda_j$ (Corollary 1). □

Then, in particular, the second largest eigenvalue of $\mathcal{P}^{(k)}$ is:

$$\lambda_2^{(k)} = \frac{N-2}{2k}\left(1 - \left(1 - \frac{2}{N}\right)^k\right) \qquad (19)$$

*Remark 12* For all $k \geq 2$, $\lambda_2 > \lambda_2^{(k)}$. First of all, we can observe that $\lambda_2^{(k)}$ decreases for increasing $k$. Therefore:

$$\lambda_2 - \lambda_2^{(k)} \geq \lambda_2 - \lambda_2^{(2)} = \frac{N-2}{4N^2}(4N - N^2 + (N-2)^2) = \frac{N-2}{N^2} > 0$$

For $N \gg 1$, the series expansion of $\lambda_2^{(k)}$ gives us:

$$\lambda_2^{(k)} = \frac{N-2}{2k}\left(1 - \left(1 - \frac{2k}{N} + \frac{2k(k-1)}{N^2} + \mathcal{O}\left(\frac{1}{N^3}\right)\right)\right)$$
$$= \frac{N-2}{N} - \frac{(N-2)(k-1)}{N^2} + \mathcal{O}\left(\frac{1}{N^2}\right)$$

We can observe how the spectral gap changes. If we consider the series expansion of $\left(1 - \frac{2}{N}\right)^k$ for $N \to \infty$, we get:

$$\lambda_1^{(k)} - \lambda_2^{(k)} = \frac{2}{N} + \frac{(N-2)(k-1)}{N^2} + \mathcal{O}\left(\frac{1}{N^2}\right)$$

It can be interesting to choose $k$ as a function of $N$. Let us consider, for example, $k = \alpha N$, with $0 < \alpha \leq 1$. In this case, we have:

$$\begin{aligned}
\lambda_2^{(\alpha N)} &= \frac{N-2}{2\alpha N}\left(1 - \left(1 - \frac{2}{N}\right)^{\alpha N}\right) \\
&\stackrel{\text{for } N \to \infty}{=} \frac{N-2}{2\alpha N}\left(1 - \left(e^{-2\alpha} + \mathcal{O}\left(\frac{1}{N}\right)\right)\right) \\
&= \frac{(N-2)\left(1 - e^{-2\alpha}\right)}{2\alpha N} + \mathcal{O}\left(\frac{1}{N}\right) \to \frac{1 - e^{-2\alpha}}{2\alpha} \text{ for } N \to \infty
\end{aligned}$$

We can observe that $\frac{1-e^{-2\alpha}}{2\alpha} =: \overline{\lambda}_2^{(\alpha N)}$ decreases when $\alpha$ increases. Moreover:

- $\overline{\lambda}_2^{(\alpha N)} \to 1$ for $\alpha \to 0$, which means that the spectral gap, $1 - \lambda_2^{(\alpha N)}$ converges to zero for $N \to \infty$ and $\alpha \to 0$;
- If $\alpha = 1$ then $\overline{\lambda}_2^{(N)} = \frac{1}{2} - \frac{1}{2e^2}$. Therefore, the spectral gap will be $\frac{1}{2} + \frac{1}{2e^2}$

The spectral gap indicates how quickly our RW converges through its stationary distribution. As expected, if $\alpha \to 0$ then the spectral gap will be close to 0. On the other hand for all $\alpha > 0$ the spectral gap tends to a strictly positive quantity, while the spectral gap corresponding to the case of the basic model converges to zero for $N \to \infty$. In particular, when $\alpha = 1$ (*i.e.* we are considering the optimal case, in which we are allowed to do among 1 and $N$ mutations at each mutation step), the spectral gap, $\frac{1}{2} + \frac{1}{2e^2}$, is significantly bigger than the one obtained for the basic model, $2/N$.

3.2 Comparison of hitting times

In this section we compare hitting times referring to some relevant mutational models we have already presented. We do not consider models that entail non-connected graphs (the model of class switch of $k$-length strings and the exchange mutation model): this choice is motivated by the discussion from the beginning of Section 3.1.3. In Table 1 below we collect most important characteristics of these RWs on $\{0,1\}^N$: the hitting time and its approximation for big $N$, that we will discuss in this current section, the stationary distribution and the value of the second larger eigenvalue when known.

Table 1: Table 1 summarizes the main characteristics of most random processes we introduce and analyze in Sections 2 and 3.

| Model | Hitting time | Stationary distribution | Second bigger eigenvalue |
|---|---|---|---|
| **Basic model** | $H(\bar{d}) = \sum_{d=0}^{\bar{d}-1} \frac{\sum_{j=1}^{N-1-d} C_N^{d+j}+1}{C_{N-1}^d} \sim 2^N$ | $\pi$ | $1 - \frac{2}{N}$ |
| **Switch 1-2** | $\sim 2^{N-1}$ | $\pi\big|_{\overline{\mathcal{H}_N}(1-2)}$ | - |
| **Allowing 1 to $k$ mutations** | $\overline{T}_N^{(k)}(\bar{d}) = \sum_{l=2}^{2^N} \mu_l^{(k)} - \frac{1}{2^N C_N^{\bar{d}}} \sum_{l=2}^{2^N} \mu_l^{(k)} R_N(l,\bar{d})$ | $\pi$ | $\frac{N-2}{2k}\left(1 - \left(\frac{N-2}{N}\right)^k\right)$ |

*3.2.1 Class switch of 1 or 2-length strings depending on the Hamming distance to $\bar{\mathbf{x}}$.*

We use results obtained in Section 2 for the $(D_n)$ process concerning the SRW on the $N$-dimensional hypercube and we apply them to this model. Here we shall introduce another definition of the distance, which is adapted to a connected component $\mathcal{H}_{N,2} \subset \{0,1\}^N$, where we denote by $\mathcal{H}_{N,2}$ one of the two parts in which $\{0,1\}^N$ is divided applying the model of class switch of 2-length strings. We recall that $\mathcal{H}_{N,2}$ is a $(N-1)$-dimensional hypercube, and that the graph underlying the model of class switch of 1 or 2-length strings corresponds essentially to the graph obtained with the model of switch of 2-length strings, up to the addition of some oriented edges from $\mathcal{H}_N^{(1-2)}$ to $\overline{\mathcal{H}_N}^{(1-2)}$.

**Definition 14** For all $\mathbf{x}_i, \mathbf{x}_j \in \mathcal{H}_{N,2}$ we denote by $h^{(2)}(\mathbf{x}_i, \mathbf{x}_j)$ the number of edges in a shortest path connecting them. Simultaneously we denote by $D_n^{(2)} = h^{(2)}(\mathbf{X}_n, \overline{\mathbf{x}})$, $D_n^{(2)} \in \{0, \ldots, N-1\}$ $\forall n \geq 0$.

Considering the process $(D_n^{(2)})_{n \geq 0}$, all the results we determined in Section 2 hold true. Furthermore, let us denote by $\mathbb{E}_{\mathbf{x}_i}^{(2)}[\tau_A]$ the expected number of steps before set $A \in \mathcal{H}_N^2$ is visited starting at $\mathbf{x}_i \in \mathcal{H}_N^2$ and following the model of switch of 2-length strings. Then, we also denote by $H_{N-1}^{(2)}(d) = \mathbb{E}_{\mathbf{x}}^{(2)}[\tau_{\{\overline{\mathbf{x}}\}}]$ where $d = h^{(2)}(\mathbf{x}, \overline{\mathbf{x}})$.

*Remark 13* Clearly if $D_0 = 2t$ and $t > 0$, which means that $\mathbf{X}_0$ and $\overline{\mathbf{x}}$ belong to the same connected component in the model of class switch of 2-length strings, then the mean hitting time for the current model will be of the order of a half the mean hitting time for the basic model, as we are considering here a $(N-1)$-dimensional hypercube instead of a $N$-dimensional one.

The result below, which is an immediate application of the Ergodic Theorem, will help us understanding better the general behavior of this mean hitting time:

**Proposition 11** *Let $(\mathbf{X}_n)_{n \geq 0}$ be a SRW on $\mathcal{H}_N$. We denote by $T_d^+ := \inf\{n \geq 1 \,|\, D_n = d\}$ and $T_d := \inf\{n \geq 0 \,|\, D_n = d\}$. Then:*

$$\mathbb{E}_{D_0 = d}[T_d^+] = \frac{2^N}{C_N^d} \quad (20)$$

*Proof* The proof is obtained by applying the Ergodic Theorem to the $(D_n)$ process and its stationary distribution, the binomial probability distribution. □

For the discussion we made in Section 2.2 and, in particular, Remark 3 we can conclude that for $N \gg 1$ the order of magnitude of the time we spend to reach the $N$ nodes at Hamming distance 1 from $\overline{\mathbf{x}}$ is:

$$\mathbb{E}_{D_0 = d}[T_1] \sim \frac{2^N}{N} \quad (21)$$

Then we can claim the following result, which comes directly from Equation (21):

**Proposition 12** *Let us suppose that $D_0 = 2t^* + 1$ with $0 < t^* \leq \lfloor (N-1)/2 \rfloor$. Then for $N \gg 1$ we have:*

$$\mathbb{E}_{D_0 = d}^{(2)}[T_1] \sim \frac{2^{N-1}}{N}$$

Finally we have:

**Proposition 13** *We denote by $\mathbb{E}_{\mathbf{x}_0}^{(1-2)}[\tau_{\{\overline{\mathbf{x}}\}}]$ the mean hitting time to reach $\overline{\mathbf{x}}$ starting from $\mathbf{x}_0$ and referring to the mutation model of class switch of 1 or 2 length strings. Then, for $N \gg 1$ we have:*

$$\mathbb{E}_{\mathbf{x}_0}^{(1-2)}[\tau_{\{\overline{\mathbf{x}}\}}] \sim \frac{1}{2}\mathbb{E}_{\mathbf{x}_0}[\tau_{\{\overline{\mathbf{x}}\}}] \quad with \quad \mathbb{E}_{\mathbf{x}_0}[\tau_{\{\overline{\mathbf{x}}\}}] \sim 2^N,$$

*where $\mathbb{E}_{\mathbf{x}_0}[\tau_{\{\overline{\mathbf{x}}\}}]$ is the hitting time from $\mathbf{x}_0$ to $\overline{\mathbf{x}}$ according to the basic model, as defined in Section 2.3.*

*Proof* First of all we observe that the last statement is a direct consequence of Proposition 3. As far as the first statement is concerned, we observe that according to the model we are analyzing here and due to Proposition 12, for $N \gg 1$ the order of magnitude of $\mathbb{E}_{\mathbf{x}_0}^{(1-2)}[\tau_{\{\overline{\mathbf{x}}\}}]$ is:

$$\mathbb{E}_{\mathbf{x}_0}^{(1-2)}[\tau_{\{\overline{\mathbf{x}}\}}] \sim \frac{1}{2}\left(\frac{2^{N-1}}{N} + 2^{N-1}\right) + \frac{1}{2}2^{N-1}$$

where the first term corresponds to the case $\mathbf{x}_0 \notin \overline{\mathcal{H}_N}^{(1-2)}$ and the second one corresponds to the opposite case (as we choose randomly the first vertex, $\mathbf{x}_0$, we have probability $1/2$ that it belongs to each part of the hypercube). For the last term we used again Proposition 3 applied to a $(N-1)$-dimensional hypercube and according to the $(D_n^{(2)})$ process and the corresponding hitting time $H_{N-1}^{(2)}(d)$. The result follows. □

Table 2: Average expected times from $[0,\ldots,0]$ to $[1,\ldots,1]$, comparing the basic mutational model and the model of class switch of 1 or 2 length strings. Here we denote by $\widehat{\tau_{\{\overline{\mathbf{x}}\}}}_n$ the average value obtained over $n$ simulations and by $\widehat{\sigma}_n$ its corresponding estimated standard deviation.

| Mutational model | $N$ | $n$ | $\widehat{\tau_{\{\overline{\mathbf{x}}\}}}_n$ | $\frac{\widehat{\sigma_n}}{\sqrt{n}}$ |
|---|---|---|---|---|
| **Basic** | 10 | 5000 | 1188.7996 | 16.2930 |
|  | 11 | 5000 | 2312.5648 | 32.1073 |
| **Switch 1-2** | 10 | 5000 | 602.8124 | 8.4773 |
|  | 11 | 5000 | 1181.5174 | 16.9023 |

*Remark 14* We simulated the basic mutational model and the model of class switch of 1 or 2 length strings in order to compare the hitting times from $\mathbf{x}_0 := [0,\ldots,0]$ to $\overline{\mathbf{x}} := [1,\ldots,1]$ for both mutational models. We consider the case $N = 10$ and $N = 11$ in order to have an example in which the process starts from $\overline{\mathcal{H}_N}^{(1-2)}$ and from $\mathcal{H}_N^{(1-2)}$ respectively. Indeed, if $N = 10$ the process starts from the connected component to which $\bar{\mathbf{x}}$ belongs, while when $N = 11$ we have to reach one of the $N$ nodes having distance 1 from $\bar{\mathbf{x}}$ to reach the connected component containing $\bar{\mathbf{x}}$. The average resulting hitting times are summarized in Table 2.

*3.2.2 Allowing 1 to k mutations.*

In this section we study the mean hitting time to cover a fixed Hamming distance $d$. First of all, we give the expression of the hitting time from node $i$ to node $j$ using the spectra. This formula is deduced by the more general one given in [36], in the case of regular graphs (the graph we obtained by a convex combination of matrices $\mathcal{P}^i$ is a regular multigraph). We refer to the notations given in Section 2 for the eigenvectors of matrix $\mathcal{P}$: $\mathbf{v}_s = (v_{s1}, \ldots, v_{s2^N})$ is the eigenvector of $\mathcal{P}$ corresponding to $\lambda_s$. These eigenvectors are the columns of matrix $Q_N$ (Section 2.1), and each component $v_{si}$ corresponds to node $i$, as they were organized while constructing the adjacency matrix. Denoting by $T(i,j)$ the hitting time from node $i$ to node $j$ in $\mathcal{H}_N^{(k)}$, we obtain the following expression:

$$T(i,j) = 2^N \sum_{l=2}^{2^N} \frac{1}{1-\lambda_l^{(k)}} (v_{lj}^2 - v_{li} v_{lj}),$$

which can be written using column vectors of $\mathcal{Z}_N$.

$$T(i,j) = \sum_{l=2}^{2^N} \frac{1}{1-\lambda_l^{(k)}} (z_{lj}^2 - z_{li} z_{lj})$$

We are interested in studying the equation below:

$$\overline{T}_N^{(k)}(d) := \frac{1}{2^N C_N^d} \sum_{h(i,j)=d} T(i,j) = \frac{1}{2^N C_N^d} \sum_{l=2}^{2^N} \frac{1}{1-\lambda_l^{(k)}} \sum_{h(i,j)=d} (z_{lj}^2 - z_{li} z_{lj}), \tag{22}$$

where $2^N C_N^d$ corresponds to the number of couples of nodes of $\{0,1\}^N$ having Hamming distance $d$.
First of all we can observe that for all $l$ and for all $j$, $z_{lj}^2 = 1$. Moreover, in order to simplify notations, we denote $\mu_l^{(k)} := (1-\lambda_l^{(k)})^{-1}$. Also, we denote $R_N(l,d) := \sum_{h(i,j)=d} z_{li} z_{lj}$. We have proved:

**Proposition 14**

$$\overline{T}_N^{(k)}(d) = \sum_{l=2}^{2^N} \mu_l^{(k)} - \frac{1}{2^N C_N^d} \sum_{l=2}^{2^N} \mu_l^{(k)} R_N(l,d) \tag{23}$$

All the elements of this equation are known, except $R_N(l,d)$. Let us consider the $2^N \times (N+1)$ matrix $\mathcal{R}_N = (R_N(l,d))$, with $1 \leq l \leq 2^N$ and $0 \leq d \leq N$. One can prove by iteration:

**Proposition 15**
$$\mathcal{R}_N = \mathcal{Z}_N \cdot \mathcal{L}_N \qquad (24)$$
*where $\mathcal{Z}_N := (\mathbf{z}_1, \ldots, \mathbf{z}_{2^N})$ is recursively obtained from $\mathcal{Z}_{N-1}$ (Section 2.1), and*

$$\begin{cases} \mathcal{L}_1 = 2I_2, \ I_n \text{ being the n-dimensional identity matrix} \\ \\ \mathcal{L}_N = \begin{pmatrix} 2 \cdot \mathcal{L}_{N-1} & \mathbf{0}_{2^{N-1}} \\ \mathbf{0}_{2^{N-1}} & 2 \cdot \mathcal{L}_{N-1} \end{pmatrix}, \ \mathbf{0}_n \text{ being the n-length zero column vector} \end{cases}$$

*3.2.3 Numerical simulations*

In Figure 3 we plot some examples of the dependence of $\overline{T}_N^{(k)}(d)$ on $d$ and $k$ for different values of $N$.

Figure 3a shows that for increasing $k$, $\overline{T}_N^{(k)}(d)$ varies on a smaller interval: $[1023, 1186.5]$ for $k=1$, $[1028.1, 1068.6]$ for $k=5$ and $[1025.6, 1044.8]$ for $k=10$. It is intuitive to understand this fact: the hitting time depends less from the initial Hamming distance if we allow to make more mutations at the same mutation step. Indeed, we can actually visit more distant nodes since the first steps, so the initial Hamming distance has a smaller influence on the result. Figures 3b and 3c show the dependence of $\overline{T}_N^{(k)}(d)$ on $k$. We obtain the best result for the biggest $k$, except in the case $d=1$ (as already shown by Figure 3a). Curves corresponding to the case $d=5$ and $d=10$ are really close: we can evaluate their minimal and maximal values, which are respectively 1043.25 and 1177.60 for $d=5$, and 1044.82 and 1186.54 for $d=10$. This fact highlights once again that if $d>1$, the initial Hamming distance poorly influences the value of the hitting time. The case $d=1$ shows surprisingly that the hitting time is not necessarily a monotone function of $k$. Figure 3c allows us to focus to this behavior and better understand its causes. Indeed, as $N$ is quite small, this figure shows more clearly the oscillating behavior of $\overline{T}_N^{(k)}(d)$ while studying its dependence on $k$: for even values of $k$, $\overline{T}_5^{(k)}(1)$ increases, while for odd values of $k$ it decreases. Intuitively, as the distance we want to cover is $d=1$, if we allow to do 2 mutations instead of simply one, then we have a high probability to go further since the beginning of the process. Let us now look to Equation (22) and, in particular to the factor: $\sum_{l=2}^{2^N}(1-\lambda_l^{(k)})^{-1}$. We can understand the phenomenon plotted in Figure 3c by looking at Proposition 10. If $k$ is odd and little enough then the last eigenvalue, which is negative (equal to $-1/k$), has an important negative influence over the value of $\overline{T}_N^{(k)}(d)$. Clearly, this fact has a substantial effect only if $N$ and $k$ are little enough, otherwise it will be compensated by the effect of all other eigenvalues.

One may wonder what would be the best choice for the coefficients $a_i$ (Definition 12), $1 \leq i \leq k$, so that $\overline{T}_N^{(k)}(d)$ is minimized for a fixed $k$. We have to minimize the convex combination $\sum_{i=1}^{k} a_i \lambda_l^i$. The answer is quite evident: if $k > 2$, then the minimum is obtained by taking all $a_i = 0$ and $a_{k^*} = 1$, where $k^* = 2\lfloor (k+1)/2 \rfloor - 1$. Then the best choice for the transition probability matrix is $\mathcal{P}^{k^*}$. The fact that we need to consider the greater odd component has also another explanation, which is more intuitive. Indeed if we consider the RW given by $\mathcal{P}^{2t}$, then due to the bipartite structure of the hypercube we will be trapped in one of the connected-components of the graph, *i.e.* the resulting graph is non-connected. That means that we will not be able to reach those nodes having a different parity of 1s in their string, referring to $\mathbf{X}_0$.

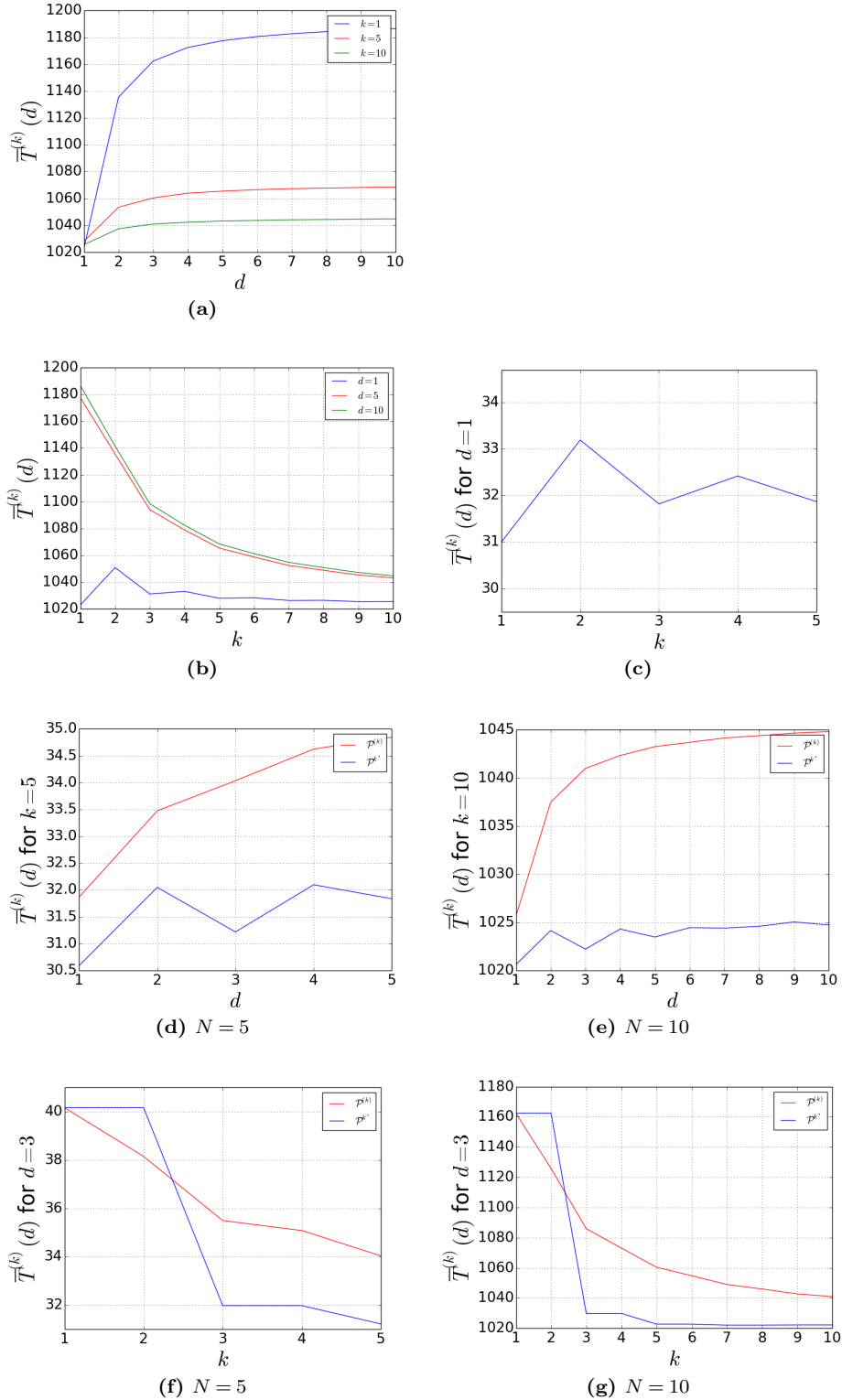

**Figure 3:** (a) Dependence of $\overline{T}_N^{(k)}(d)$ on $d$ for $N=10$ and $k=1$, 5 or 10. (b) Dependence of $\overline{T}_N^{(k)}(d)$ on $k$ for $N=10$ and different values of $d$. (c) Dependence of $\overline{T}_5^{(k)}(1)$ on $k$. (d, e) Dependence of $\overline{T}_N^{(k)}(d)$ on $d$ for different values of both $N$ and $k$. Values obtained by using as transition probability matrices $\mathcal{P}^{(k)}$ and $\mathcal{P}^{k^*}$ respectively are compared. (f, g) Dependence of $\overline{T}_N^{(k)}(d)$ on $k$ for different values of both $N$ and $d$. Again, cases corresponding to $\mathcal{P}^{(k)}$ and $\mathcal{P}^{k^*}$ are compared.

In Figures 3d, 3e, 3f and 3g we plotted together the values of the hitting time to cover a Hamming distance $d$ for different values of $N$, $k$, and $d$, comparing the process given by $\mathcal{P}^{(k)}$ and the one corresponding to $\mathcal{P}^{k^*}$. This gives more evidence of the fact that the second one is the optimal one. It is interesting to look at the case in which $d$ is fixed and we let $k$ vary. For $k = 1$ both processes gave the same result as $\mathcal{P}^{1^*} = \mathcal{P} = \mathcal{P}^{(1)}$. Moreover we necessarily have that for $k = 2$ the process $\mathcal{P}^{(2)}$ is the faster one: we recall that defining $\mathcal{P}^{k^*}$ we consider the greater odd $k$, and then $\mathcal{P}^{2^*} = \mathcal{P}$, while the process $\mathcal{P}^{(2)}$ allows to do 1 or 2 mutations at each mutation step. Then $\mathcal{P}^{k^*}$ is actually the best choice among all possible convex combinations of $\mathcal{P}^i$ iff $k > 2$. In Figures 3d and 3e we observe the oscillating behavior of $\overline{T}_N^{k^*}(d)$. That depends on the structure of $\mathcal{R}_N$, considering that $\sum_{l=2}^{2^N-1} R_N(l,d) = 0$ for $d$ odd and $\sum_{l=2}^{2^N-1} R_N(l,d) = -2(2^N C_N^d)$ for $d$ even. One can get convinced of this fact by explicitly compute $\overline{T}_N^{k^*}(d)$ for $N = 3$. Moreover simulations show that this behavior is softened for increasing $d$, and that $\overline{T}_N^{k^*}(N-1) > \overline{T}_N^{k^*}(N)$. This fact is actually confirmed by simulations on the real process. Finally, Figures 3f and 3g clearly show that for $k = 2$ the process given by $\mathcal{P}^{(k)}$ allows to cover quickly a fixed Hamming distance. As expected, the best hitting time is obtained for $k = N$, and for increasing $N$ and $k$ the value of this hitting time has a smaller variation.

Table 3: An example of comparison between the theoretical and experimental values of $\overline{T}_5^{(5)}(4)$ for $\mathcal{P}^{(5)}$. $\widehat{\overline{T}_5^{(5)}(4)}_n$ denotes the average value obtained over $n$ simulations and $\widehat{\sigma}_n$ its corresponding estimated standard deviation.

| Transition probability matrix | $N$ | $d$ | $k$ | $n$ | $\overline{T}_5^{(5)}(4)$ | $\widehat{\overline{T}_5^{(5)}(4)}_n$ | $\frac{\widehat{\sigma}_n}{\sqrt{n}}$ |
|---|---|---|---|---|---|---|---|
| $\mathcal{P}^{(k)}$ | 5 | 4 | 5 | 480000 | 34.62 | 34.67 | 0.05 |

We can test all these observations by simulating the real process for both transition probability matrices, $\mathcal{P}^{k^*}$ and $\mathcal{P}^{(k)}$. Results obtained are consistent with our theoretical analysis. In order to give an idea of the values we can obtain by testing the process, in Table 3 we compare the theoretical value of $\overline{T}_N^{(k)}(d)$ corresponding to $\mathcal{P}^{(k)}$, and the experimental value with its precision, for $N = 5$, $k = 5$ and $d = 4$.

## 4 Modeling issues

The mathematical framework described in previous sections can be used to model mutations characteristic of SHM. In Sections 4.1 and 4.2 we give some more details about GCs and the binding between B-cells and antigens. Therefore, in Section 4.3 we set the modeling assumptions which justify to mathematically describe SHM process as RW on binary strings. Of course, this is a not exhaustive approximation. Hence, some limitations are discussed in Section 4.4 and some propositions for further developments are given as well.

### 4.1 The Germinal center reaction

Antigen-activated B-cells, together with their associated T cells, move into a primary lymphoid follicle, where they proliferate and ultimately form a GC. GCs are composed mainly of B-cells, but antigen specific T-cells, which have also been activated and migrated to the lymphoid follicle, make up about 10% of GC lymphocytes and provide indispensable help to B-cells. Indeed, when B-cells start to proliferate in GC, they need to receive proper survival signals, or they die by apoptosis. The number of B-cells within a germinal center grows at high pace: it can double every 6-8 hours. After about 3 days of strong proliferation, B-cells start undergoing SHM, in order to diversify the variable region of their BCRs, and those cells that express newly generated modified BCRs are selected for enhanced antigen binding. The fast proliferation rate of B-cells is required for the generation of a large number of modified BCRs within a short frame time (one cell gives $10^4$ blasts in 72 hours). Some B-cells positively selected in the light zone differentiate into memory B-cells or plasma cells. The GC reaches its maximal size within approximately two weeks, after which the structure slowly involutes and disappears within several weeks. During the GC process B-cells are subjected to powerful selection mechanisms that facilitate the generation of high affinity antibodies: a B-cell that express a newly generated BCR needs to be tested for enhanced antigen binding. This process is mediated by FDCs and follicular helper T-cells. BCR stimulation through antigen binding coupled with co-stimulatory signals that are transmitted to the B-cell by GC T-cells, provides survival signals to the cell; by contrast, failure of the BCR to bind antigen causes cell death by apoptosis. The final differentiation of a GC B-cell into a plasma cell or a long-lived memory B-cell is driven by the acquisition of a high-affinity BCR. For short-lived memory B-cells, the differentiation process seems to be stochastic, as throughout GC formation GC B-cells are constantly selected to enter the memory pool [40, 52].

4.2 B-cell receptors and antigen-antibody binding

Immunoglobulins (Ig) present at the antigen receptor are Y-shaped macro proteins composed of four polypeptide chains assembled by disulfide bonds: two identical heavy (H) chains and two identical light (L) chains. Each chain consists of two regions: a constant (C) region, which has an effector function, and a variable (V) region composed by the variable parts of the two chains together. During GC reaction the only one involved in SHMs is the V region, which also determines the antigen binding site. We call *antigen binding site* or *paratope* the specialized portion of the BCR V region used for identifying other molecules, while the regions on any molecule that the paratopes can recognize are called *epitopes*. B-cells are able to bind ligands whose surfaces are 'complementary' to that of their antigen binding site, where complementarity means that the amino-acids composing the paratope and the epitope are distributed in such a way to form bonds which are able to hold the antigen to the B-cell. In this case these bonds are all non-covalent (as hydrogen bonds, electrostatic bonds, van der Waals forces and hydrophobic bonds), which are by their nature reversible. Multiple bonding between the antigen and the B-cell ensures that the antigen is bound tightly to the B-cell. The interaction between paratope and epitope can be characterized in terms of a binding affinity, that will be proportional to their complementarity. The *affinity* is the strength of the reaction between a single antigenic determinant and a single combining site on the B-cell: it summarizes the attractive and repulsive forces operating between the antigenic determinant and the combining site of the B-cell, and corresponds to the equilibrium constant that describes the antigen-B-cell reaction [19].

Each antigen typically has several epitopes, so that the surface of an antigen presents variable motifs that B-cells, through their receptors, can discriminate as distinct epitopes. If we define an epitope by its spatial contact with a BCR during binding, the number of relevant amino-acids is approximately 15, and among these amino-acids only around 5 in each epitope strongly influence the binding. These strong sites may contribute about one-half of the total free energy of the reaction, while the other amino-acids influence in binding constant by up to one order of magnitude or even have no detectable effect.

Simultaneously, a BCR contains a variety of possible binding sites and each antibody binding site defines a paratope: about 50 variable amino-acids make up the potential binding area of a BCR. In agreement with the above, only around 15 among these 50 amino-acids physically contact a particular epitope: these define the structural paratope. Consequently, antibodies have a large number of potential paratopes as the 50 or so variable amino-acids composing the binding region define many putative groups of 15 amino-acids.

Substitutions both in and away from the binding site can change the spatial conformation of the binding region and affect the binding reaction. The consequence of mutation at a particular site depends on the original amino-acid and the amino-acid used for substitution ([19], Chapter 4).

4.3 From DNA to amino-acids: choosing the best viewpoint

Mutations observed on the binding site of B-cells during the GC process are the result of genetic mutations produced by SHM on the portion of DNA encoding for the BCR V region. In the current section we discuss a model of genetic mutations and its effects on the amino-acid string, under the assumption of having two amino-acid classes. We show that the framework we set up in previous sections can adapt to model the effects of SHM over BCRs and study the variation of the affinity with the presented antigen.

The genetic code is a sequence of four nucleotides, guanine (G), adenine (A) (called purines), thymine (T) and cytosine (C) (pyrimidines), joined together. They make three-letter words: the codons. Each codon corresponds to a specific amino-acid or to a stop signal, which interrupts the building of the protein during translation. As the number of possible combinations of 4 nucleotides in 3-length words is 64, and there exists 20 amino-acids in naturally derived proteins, more than a single codon codes for the same amino-acid [51]. Table 4 shows the correspondence between codons and amino-acids.

Table 4: The correlation between codons and amino-acids: most of the amino-acids derives from more than a single codon.

|   |     | **T** |     |     | **C** |     |     | **A** |     |     | **G** |     |     |
|---|-----|-------|-----|-----|-------|-----|-----|-------|-----|-----|-------|-----|-----|
| **T** | TTT | Phe (F) | TCT | Ser (S) | TAT | Tyr (Y) | TGT | Cys (C) | **T** |
|   | TTC | Phe (F) | TCC | Ser (S) | TAC | Tyr (Y) | TGC | Cys (C) | **C** |
|   | TTA | Leu (L) | TCA | Ser (S) | TAA | Stop | TGA | Stop | **A** |
|   | TTG | Leu (L) | TCG | Ser (S) | TAG | Stop | TGG | Trp (W) | **G** |
| **C** | CTT | Leu (L) | CCT | Pro (P) | CAT | His (H) | CGT | Arg (R) | **T** |
|   | CTC | Leu (L) | CCC | Pro (P) | CAC | His (H) | CGC | Arg (R) | **C** |
|   | CTA | Leu (L) | CCA | Pro (P) | CAA | Gln (Q) | CGA | Arg (R) | **A** |
|   | CTG | Leu (L) | CCG | Pro (P) | CAG | Gln (Q) | CGG | Arg (R) | **G** |
| **A** | ATT | Ile (I) | ACT | Thr (T) | AAT | Asn (N) | AGT | Ser (S) | **T** |
|   | ATC | Ile (I) | ACC | Thr (T) | AAC | Asn (N) | AGC | Ser (S) | **C** |
|   | ATA | Ile (I) | ACA | Thr (T) | AAA | Lys (K) | AGA | Arg (R) | **A** |
|   | ATG | Met (M) | ACG | Thr (T) | AAG | Lys (K) | AGG | Arg (R) | **G** |
| **G** | GTT | Val (V) | GCT | Ala (A) | GAT | Asp (D) | GGT | Gly (G) | **T** |
|   | GTC | Val (V) | GCC | Ala (A) | GAC | Asp (D) | GGC | Gly (G) | **C** |
|   | GTA | Val (V) | GCA | Ala (A) | GAA | Glu (E) | GGA | Gly (G) | **A** |
|   | GTG | Val (V) | GCG | Ala (A) | GAG | Glu (E) | GGG | Gly (G) | **G** |

Different kind of genetic mutations can affect the DNA sequence of a gene. They can be regrouped in three main categories: base substitutions, inser-tion and deletions. A single base substitution is a switch of a nucleotide with another. This is the simplest kind of mutation and it can turn out to be missense, nonsense or silent, once we observe the resulting new protein. We said that a mutation is missense if the result of the genetic mutation is a different amino-acid in the protein. The mutation is nonsense when the genetic mutation results in a stop codon instead of an amino-acid. Finally, a silent mutation is a mutation with no effect on the amino-acid string, *i.e.* the mutated sequence codes for an amino-acid with identical binding properties. We talk about insertion (resp. deletion) when one or more nucleotides are added (resp. removed) at some place within the DNA code. These last kinds of mutations can both be frameshift mutations, which are given by the insertion or deletion of a number of bases that is not a multiple of 3, altering the reading frame of the gene. SHM introduces mostly single nucleotide exchanges, together with small deletions and duplications, *i.e.* the insertion of extra copies of a portion of genetic material already present within the DNA code [23, 11]. Among these point mutations, transitions (*i.e.* substitution of a purine nucleotide with another purine one, or a pyrimidine with a pyrimidine) dominate over transversions (substitution of a purine with a pyrimidine or conversely). About half of the mutations (53%) have been estimated to be silent, about 28% nonsense, and only about 19% of all mutations have been estimated to be missense and then have an effect on affinity, which can either be of an improving nature, or of worsening and even lead to the formation of autoreactive clones [49, 26, 35].

The 20 existing amino-acids are typically classified in charged amino-acids, polar (non-charged) amino-acids and hydrophobic amino-acids, depending on their chemical characteristics. As we already discussed in Section 4.2 the bonding between BCR and antigen is made thanks to non-covalent bonding, in particular ionic bonds and hydrogen bonds. Ionic bonds are the result of the interactions between two amino-acids oppositely charged: arginine (R) and lysine (K) are positively charged, while aspartic acid (D) and glutamic acid (E) are negatively charged. As long as hydrogen bonds are concerned, also polar amino-acids can participate. In particular arginine (R), lysine (K) and tryptophan (W) have hydrogen donor atoms in their side chains; aspartic acid (D) and glutamic acid (E) have hydrogen acceptor atoms in their side chain while asparagine (N), glutamine (Q), histidine (H), serine (S), threonine (T) and tyrosine (Y) have both hydrogen donor and acceptor atoms in their side chains.

Stop codons also have an important role. Indeed, during translation (the last step necessary to build a protein starting from the DNA molecule) amino-acids continue to be added until a stop codon is reached. There exists two types of mutations involving stop codons, named nonsense and nonstop respectively. The first one corresponds to the substitution of an amino-acid with a stop codon, while the second one is the opposite case. In both cases the resulting protein has an abnormal length, which often causes a loss of function. Moreover, errors given by both nonsense and nonstop mutations are linked to over 10% of human genetic diseases [8].

Concerning mutation in activated B-cells, SHM is driven by an enzyme called activation-induced cytidine deaminase (AID) which is expressed specifically in this case. This protein can bind to single-stranded DNA only. Thus it seems to target only genes being transcribed (for which the transcription phenomenon separates temporarily double stranded DNA into small portions of two single stranded DNA sequences) [28]. AID converts Cytosine (C) in Uracil (U) by deamination. This substitution occurs at higher rates in hot spots motives like $D\underline{G}YW/WR\underline{C}H$ where ($G:C$ is the mutable position and $D \in \{A,G,T\}$, $H \in \{A,C,T\}$, $R \in \{A,G\}$, $W \in \{A,T\}$ and $Y \in \{C,T\}$, and the underlined letters are the loci of mutations) [47]. Then, two mechanisms tend to repair lesions in the DNA caused by these substitutions of C by U [10] :

a) either *mismatch repair* : substitution for the damaged zone by another sequence of nucleotides thanks to proteins MSH 2/6. The $U$ base is read as $T$ leading to a transition from a $C:G$ pair to $T:A$.
b) or *base excision repair* : U is excised by a successive actions of uracil-DNA glycolase (UNG) and apurinic/apyrimidinic endonuclease (APE1). The DNA contains then a nick, after replication of a random nucleotide is inserted in order to fill the vacant space leading to transversions and transitions.

From the mathematical point of view this is equivalent to define the switch with a random nucleotide depending on the motives present in the chain. The probability concerning the choice of this nucleotide to be inserted shall not be uniform due to the presence of mismatch and excision repairs [10]. This is not taken into account in the model we developed.

We can therefore make the following three main assumptions to model the SHM process acting on the BCR V region:

*Modeling assumption 1* SHM introduces only single point mutations in the DNA strand, missense or silent. Therefore we do not take into account nonsense mutations, in order to avoid an interruption of the mutation process due to the introduction of a stop codon. The choice of the base used for substitution is made randomly, without considering that we have mostly $A \leftrightarrow T$ and $T \leftrightarrow C$ substitutions.

*Modeling assumption 2* We consider only electrostatic and hydrogen bonds as responsible for the bonding between BCR and antigen. We suppose we have two amino-acid classes represented as 0 and 1 respectively: we denote by 1 those amino-acids which have hydrogen donor atoms in their side chains (or which are positively charged) and by 0 those amino-acids which have hydrogen acceptor atoms in their side chains (or which are negatively charged). We arbitrary chose to assign 0 or 1 to amino-acids which can act as an acid or a base in hydrogen bonds. As an exemple, as serine can form hydrogen bonds with arginine and threonine, one can assign 0 to serine and 1 to threonine (arginine is represented by 1 as it is positively charged). While translating the amino-acid chain into a binary chain, we omit all hydrophobic amino-acids, as they do not participate in electrostatic or hydrogen bonds. Their position corresponds to an empty case, which does not contribute to the affinity between B-cell and antigen. This is clearly an important simplification we make in order to build this mathematical model. We will further discuss this choice in Section 4.4.

*Modeling assumption 3* We consider a linear contact between two amino-acid strings, without taking into account the geometrical configuration of both the BCR and the antigen.

The process starts from a DNA chain coding for a BCR, $\mathbf{X}_0^{dna}$; from which we can obtain the corresponding amino-acid chain, $\mathbf{X}_0^{aa}$ (Table 4) and, consequently, its binary expression, $\mathbf{X}_0^{bin}$.

*Exemple 1*

– $\mathbf{X}_0^{dna}=$ (GTT, GAG, CTA, GTG, GAA, AGT, GGA, GCC, GAA, GTA, AAA, AAG, CCA, GGT, AGT, AGT, GTT, AAA, GTC, AGT, TGT, AAA, GCA)

– $\mathbf{X}_0^{aa}=$ (V, Q, L, V, E, S, G, A, E, V, K, K, P, G, S, S, V, K, V, S, C, K, A)

– $\mathbf{X}_0^{bin}=(-,1,-,-,0,0,-,-,0,-,1,1,-,-,0,0,-,1,-,0,0,1,-)$

*Notation 1* Given a vector $\mathbf{X}$, we denote by $|\mathbf{X}|$ its length (counting also the empty cases, if there are some). Equivalently, given a set $\mathcal{S}$, we denote by $|\mathcal{S}|$ its size

We can formalize the translation of the nucleotides chain into the amino-acids chain as follows.

**Definition 15** Let $\mathcal{N}$ and $\mathcal{A}$ be two sets of letters with size respectively $|\mathcal{N}| = k_1$ and $|\mathcal{A}| = k_2$. Let $l$ be an integer positive number so that $k_1^l \geq k_2$. Then we define $f_{k_1,k_2,l} : \mathcal{N}^l \to \mathcal{A}$, which associate at least an $l$-length sequence of letters belonging to $\mathcal{N}$ to a letter in $\mathcal{A}$.

In our specific case, following definition 15, $\overline{\mathcal{N}} := \{$G, A, T, C$\}$ is the set of nucleotides, while $\overline{\mathcal{A}}$ is the set containing all possible amino-acids, together with the stop signal. Therefore $\overline{k_1} = 4$ and $\overline{k_2} = 21$. Moreover we know that $\overline{l} = 3$ and the function $\overline{f}_{4,21,3}$ is detailed in Table 4.

*Remark 15* We can easily observe that $\overline{l} = \min\left\{n \in \mathbb{N} | \overline{k_1}^n \geq \overline{k_2}\right\}$. Indeed, having 4 nucleotides available to build a DNA strand, we need to read them at least by 3-length blocks in order to be able to synthesize all 20 amino-acids. Moreover, choosing this value for the parameter $l$ avoids to have too many sequences of nucleotides coding for the same amino-acid.

At the beginning of the process also the antigen string in its three representations is given: $\overline{\mathbf{x}}^{dna}$, $\overline{\mathbf{x}}^{aa}$ and $\overline{\mathbf{x}}^{bin}$, with $|\mathbf{X}^{dna}| = |\overline{\mathbf{x}}^{dna}| =: 3N$. Antigen representing strings remain unchanged. At each time step a single point mutation (missense or silent) is introduced in the DNA chain coding for the BCR. So, if $\mathbf{X}_t^{dna}$ is the DNA code at time $t$, we randomly choose an index $i \in \{1, \ldots, 3N\}$, a letter $a \in \overline{\mathcal{N}}$ and we place $(X_{t+1}^{dna})_i := a$. If the new codon is a stop codon, then we choose $a' \in \overline{\mathcal{N}} \setminus \{a\}$ and we put $(X_{t+1}^{dna})_i := a'$, and so on.

In order to test the affinity, we consider the binary expression of both the BCR and the antigen, which we take in its complementary form, *i.e.* $\overline{\mathbf{x}}'^{bin} := (1 - \overline{x}_1^{bin}, \ldots, 1 - \overline{x}_{|\overline{\mathbf{x}}^{bin}|}^{bin})$. This leads us back to the definition of affinity we made in Section 2: 0 matches with 0 and 1 with 1.

Assumptions 1-3 imply that for all $t \geq 0$, $|\mathbf{X}_t^{bin}| = |\overline{\mathbf{x}}^{bin}| = N$. As we consider a linear contact between $\mathbf{X}_t^{bin}$ and $\overline{\mathbf{x}}'^{bin}$, at the positions where either $\mathbf{X}_t^{bin}$ or $\overline{\mathbf{x}}'^{bin}$ has an hydrophobic amino-acid, we suppose that no match is possible. Therefore we can extend Definition 4 of the Hamming distance in a very natural way to this more general case:

**Definition 16** We denote by $Hy(\mathbf{X}_t^{bin})$ (resp. $Hy(\overline{\mathbf{x}}'^{bin})$) the set of the indices corresponding to hydrophobic amino-acids in $\mathbf{X}_t^{bin}$ (resp. in $\overline{\mathbf{x}}'^{bin}$). Therefore the Hamming distance between $\mathbf{X}_t^{bin}$ and $\overline{\mathbf{x}}'^{bin}$ is given by:

$$h(\mathbf{X}_t^{bin}, \overline{\mathbf{x}}'^{bin}) = \sum_{\substack{i \in \{1, \ldots, N\} \\ i \notin Hy(\mathbf{X}_t^{bin}) \cup Hy(\overline{\mathbf{x}}'^{bin})}} \delta_i + |Hy(\mathbf{X}_t^{bin}) \cup Hy(\overline{\mathbf{x}}'^{bin})|$$

$$\text{where } \delta_i = \begin{cases} 1 & \text{if } (X_t^{bin})_i \neq (\overline{x}'^{bin})_i \\ 0 & \text{otherwise} \end{cases}$$

Then, for all $t \geq 0$:

$$|Hy(\mathbf{X}_t^{bin}) \cup Hy(\overline{\mathbf{x}}'^{bin})| \leq h\left(\mathbf{X}_t^{bin}, \overline{\mathbf{x}}'^{bin}\right) \leq N$$

We consider that the optimal clone is reached when:

$$aff\left(\mathbf{X}_t^{bin}, \overline{\mathbf{x}}'^{bin}\right) := N - |Hy(\overline{\mathbf{x}}'^{bin})|$$

The effects of nucleotides exchanges on the binary expression of BCRs can be multiple:

**No detectable effect** : this is the result of either a silent mutation or a missense mutation which substitutes an amino-acid with another one belonging to the same amino-acid class.

**Class-switch** , derived from a missense mutation which leads to the substitution of an amino-acid with another one belonging to the other amino-acid class.

We can further complexify this model by replacing Assumption 1 with the following one:

*Modeling assumption 4* SHM introduces mostly single point mutations in the DNA, missense or silent. With weak probability, deletions or insertions can occur. For the sake of simplicity, we suppose that a deletion (resp. an insertion) consist in the elimination (resp. the addition) of a non-stop codon. Moreover, in order to avoid the problem of a variation of the length of BCR representing strings, when a deletion occur, those bits situated on the right of the deleted one shift to the left, and a random extra codon is added at the right bottom. Conversely, if an insertion occurs, the right bottom bit is deleted.

Even if these mutation events are rare, they have remarkable effects over the structure of the underlying graph. Indeed a deletion or an insertion entails a great jump in the affinity function by producing a shift of a portion of the BCR representing string. This is not the case if we consider only single point mutations. Therefore, under Assumption 4 the graph we obtain is much more complex and allows long range random connections.

### 4.3.1 Numerical simulations

In order to evaluate how deletions and insertions affect the mean number of mutation steps to reach the desired B-cell trait, we make some numerical simulations. We compare a model in which only single point mutations are allowed to another one in which also deletions and insertions can occur. We refer to Assumption 4 to define these mutational events.

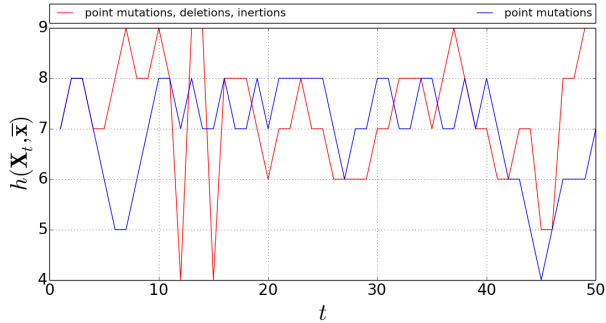

**Figure 4:** Variation of the Hamming distance through $\overline{\mathbf{x}}'^{bin}$, comparing the model of single point mutations to the one which includes also deletions and insertions (50% of all mutation events). In both cases $N = 10$. Deletions and insertions lead to a quick change in the Hamming distance. Between time 10 and 20, we can observe four deletions or insertions.

Figure 4 shows the effects of deletions and insertions over the affinity. In order to do these simulations, we arbitrary fixe a BCR and an antigen with given affinity. We do not consider those base substitutions leading to no detectable effect, *i.e.* at each time step we can observe a variation of the affinity function. In Figure 4 we can clearly locate at what time an insertion or a deletion has occurred, because this coincides with a jump of the Hamming distance between BCR and antigen.

One can ask how these random long range connections affect the average time to reach the antigen target string. Simulations show that one needs a more long time to reach $\overline{\mathbf{x}}'^{bin}$ if the probability of making such mutations increases. The results obtained through 10000 simulations are collected in Table 5.

Table 5: Average number of mutations needed to reach $\overline{\mathbf{x}}'^{bin}$, for $N = 10$ and starting from a Hamming distance 7. In $\overline{\mathbf{x}}'^{bin}$, only 2 amino-acids are hydrophobic, so by Definition 16, the optimal affinity one can reach is 8. We compare three models: in the first one no deletions nor insertions are allowed. In the second model 10% of all mutations are deletions or insertions, 50% in the last one. We denote by $\widehat{\tau_{\{\overline{\mathbf{x}}'^{bin}\}_n}}$ the average value obtained over $n$ simulations and by $\widehat{\sigma}_n$ its corresponding estimated standard deviation. Simulations show that $\widehat{\tau_{\{\overline{\mathbf{x}}'^{bin}\}_n}}$ increases when the pourcentage of deletions or insertions grows, and so does the corresponding variation.

| % deletions/insertions | $|\overline{\mathbf{x}}'^{bin}|$ | $h(\mathbf{X}_0^{bin}, \overline{\mathbf{x}}'^{bin})$ | $n$ | $\widehat{\tau_{\{\overline{\mathbf{x}}'^{bin}\}_n}}$ | $\widehat{\frac{\sigma_n}{\sqrt{n}}}$ |
|---|---|---|---|---|---|
| **0** | 10 | 7 | 10000 | 8824.93 | 86.80 |
| **10** | 10 | 7 | 10000 | 9091.12 | 92.01 |
| **50** | 10 | 7 | 10000 | 10075.89 | 100.59 |

We can discuss which viewpoint is the most suitable to study mutations and their effects over the interactions between BCR and antigen. It is really hard to define a clear correspondence between genetic mutations and the evolution of the affinity, even while considering a simple linear contact between molecules (hence without observing the changes in the geometrical structure of the protein). Indeed, in order to test the affinity between BCR and antigen we constantly need to project the DNA string on the smaller state-space containing the binary representations of B-cell traits. If we directly consider mutations on the binary strings, then the process we obtain is faster, as we do not observe missense mutations, and the evaluation of the affinity is immediate.

The comprehension of the nature of genetic mutations and their consequences on the new generated protein, suggested us to make Assumptions 1-3 to formalize the model. In particular, we found reasonable to look directly to the amino-acid chains and their binary representation, which allows to study the affinity between BCR and antigen using the Hamming distance. Therefore, under these hypotheses the general mathematical framework described in Section 2 can be applied to study how different kind of missense mutations affects the dynamics of AAM. As we show in Sections 2-3, this already brings interesting and complexes mathematical problems.

4.4 Limitations and extensions

In this paper we propose and study mutational processes on $N$-length binary strings, which can be variously applied to evolutionary contexts. As far as the application to the SHM process is concerned, we can make some remarks about our assumptions, which can bring us to enrich and complexify the model through a more coherent representation of the true biological process.

First of all we decided to consider only two amino-acid classes. From one side this assumption is justified as charged and polar amino-acids are effectively the most responsible in creating bonds that determine the antigen-antibody interaction. Therefore they strongly influence the affinity between BCR and antigen. Nevertheless, by making this simplification we omit all hydrophobic amino-acids from the string, and that is not without consequences. The elimination of hydrophobic amino-acids from the string significantly changes the structure of the chain, therefore the ability for charged and polar amino-acids to be in contact with each-others. Moreover, the effects of genetic mutations on the new generated protein could be even more complex than the ones we considered in this paper. Finally, by taking into account also hydrophobic amino-acids, we would be able to consider hydrophobic bonds, which also influences the antigen-antibody interaction. Therefore it seems more appropriate to consider three amino-acids classes, and define an affinity function so that positively charged amino-acids match with negatively charged, and hydrophobic amino-acids match with hydrophobics.

As far as the nature of mutations is concerned, we essentially described mutational processes given by combinations of single point mutations mechanisms. During SHM nucleotide exchanges are the most frequent among all possible mutations. Despite this, also some deletions and insertions occur. This has two main consequences. Firstly that means that the length of the BCR representing string could actually change during the process, while we consider it as fixed and equal to the length of the antigen. We can maybe overcome this problem by saying that the chain represented in our model corresponds to the portion of BCR in contact with the antigen, and this is almost fixed (Section 4.2). Moreover these mutations can imply substantial changes into the amino-acid chain, hence they can bring a great jump of the affinity to the presented antigen. Therefore, even if these are rare mutational events, they may have an important effect in AAM and consequently it could be interesting to takes also insertions and deletions into account. All these observations lead interesting mathematical questions. Of course we can also envisage developments in other directions. For example by considering the creation of bonds among amino-acids of the BCR (resp. the antigen) itself, which determines the geometrical structure of the protein and consequently the portion of the BCR and the antigen that can actually be in contact.

We propose some numerical simulations to evaluate the consequences over the hitting time of both the addiction of extra amino-acids classes and the possibility of having a BCR string longer than the antigen one.

A. S. Perelson and G. Weisbuch in [44] proposed a model with 3 amino-acid classes: hydrophobic, hydrophilic positively charged and hydrophilic negatively charged. Hydrophobic amino-acids match with hydrophobic and hydrophilic positively charged with hydrophilic negatively charged. We simulated the expected time to reach a given configuration comparing the model with 2 amino-acid classes and the one with 3 amino-acid classes, and considering single switch-type mutations. We take two random 10-length strings having maximal distance between each-others. We extended Definition 4 of Hamming distance to the state-space $\{0,1,2\}^N$ in a natural way, keeping the same notation: $\forall \mathbf{x} = (x_1,\ldots,x_N), \mathbf{y} = (y_1,\ldots,y_N) \in \{0,1,2\}^N$, their Hamming distance is given by:

$$h(\mathbf{x},\mathbf{y}) = \sum_{i=1}^{N} \delta_i \qquad \text{where} \qquad \delta_i = \begin{cases} 1 \text{ if } & x_i \neq y_i \\ 0 \text{ otherwise} \end{cases} \qquad (25)$$

Therefore the affinity is defined as in Definition 3. We simulated for both cases a single switch-type mutational model (Definition 5 for 2 amino-acid classes and Definition 17 below for 3 amino-acid classes), testing the time we need to reach the target vertex.

**Definition 17** Let $\mathbf{X}_n \in \{0,1,2\}^N$ be the BCR at step $n$. Let $i \in \{1,\ldots,N\}$ be a randomly chosen index, and $a \in \{0,1,2\} \setminus \{X_{n,i}\}$ a randomly chosen number. Then $\mathbf{X}_{n+1} := (X_{n,1},\ldots,X_{n,i-1},a,X_{n,i+1},\ldots,X_{n,N})$.

Table 6 shows the results we obtained over 10000 simulations.

Table 6: Average expected times to cover a Hamming distance $h(\mathbf{X}_0,\overline{\mathbf{x}}) = 10 = N$, comparing the model with 2 amino-acid classes and the one with 3 amino-acid classes. Here we denote by $\widehat{\tau_{\{\overline{\mathbf{x}}\}}}_n$ the average value obtained over $n$ simulations and by $\widehat{\sigma}_n$ its corresponding estimated standard deviation.

| Amino-acid classes | $N$ | $h(\mathbf{X}_0,\overline{\mathbf{x}})$ | $n$ | $\widehat{\tau_{\{\overline{\mathbf{x}}\}}}_n$ | $\frac{\widehat{\sigma_n}}{\sqrt{n}}$ |
|---|---|---|---|---|---|
| **2** | 10 | 10 | 10000 | 1213.2108 | 12.0138 |
| **3** | 10 | 10 | 10000 | 62160.8263 | 635.0458 |

We already knew from theoretical analysis that the order of magnitude for the hitting time of the basic mutational model is $2^N$ for $N$ big enough. Simulations clearly show that when we consider 3 amino-acid classes, the order of magnitude of the hitting time of a single switch-type mutational model significantly increases, and is of the order of $3^N$, as proved by Proposition 4. Moreover we observe that the variance corresponding to the second model is significantly bigger as well.

It is clear that if we consider more amino-acid classes, it takes much longer to reach a precise element of the new state-space. Nevertheless, one can understand that if we keep the same distance function as defined in Equation (25), than we are actually asking for a higher degree of precision while building the B-cell trait. Therefore, we can not directly compare hitting times corresponding to a model with a greater number of amino-acid classes and keeping the same affinity function as the one used with only two amino-acid classes. If one want to obtain a comparable result by using more than two amino-acid classes, one has to use a weaker definition of affinity.

**Definition 18** Let $\mathcal{S}$ be a set of letters, $|\mathcal{S}| = s > 2$. Let us partition $\mathcal{S}$ into two subsets: $\mathcal{S} := \mathcal{S}_1 \sqcup \mathcal{S}_2$. $\forall\, \mathbf{x},\, \mathbf{y} \in \mathcal{S}^N$, their distance is given by:

$$h_{\mathcal{S}_1,\mathcal{S}_2}(\mathbf{x},\mathbf{y}) = \sum_{i=1}^{N} \delta_i \quad \text{where} \quad \delta_i = \begin{cases} 1 \text{ if } x_i \in \mathcal{S}_1,\, y_i \in \mathcal{S}_2 \text{ or conversely} \\ 0 \text{ otherwise} \end{cases}$$

Consequently, their affinity is given by:

$$aff(\mathbf{x},\mathbf{y}) = N - h_{\mathcal{S}_1,\mathcal{S}_2}(\mathbf{x},\mathbf{y})$$

By using this new affinity function we can actually compare the hitting times and the order of magnitude is clearly the same.

Let us now go back to Assumption 2 and to the structure of the string given in Section 4.3 (in particular, hydrophobic amino-acids are represented by empty cases). Contrary to what stated in Assumption 4, we suppose that the BCR length can be modified by insertions and deletions. Consequently, also a modification of the distance function is needed. We arbitrary fix a BCR and an antigen with given affinity. We do not consider those base substitutions leading to no detectable effect, *i.e.* at each time step we can observe a variation of the affinity function. We suppose that 90% of all mutation events are single point mutations, 10% deletions or insertions. If we are in this case and $|\mathbf{X}_t^{bin}| > |\overline{\mathbf{x}}'^{bin}|$, then with probability $1/2$ a deletion occurs and with probability $1/2$ an insertion occur. Otherwise, it will be necessarily an insertion (this is to avoid to obtain $|\mathbf{X}_t^{bin}| = 0$). As long as the affinity is concerned, if $|\mathbf{X}_t^{bin}| > |\overline{\mathbf{x}}'^{bin}|$, $|\mathbf{X}_t^{bin}| := n_1$, $|\overline{\mathbf{x}}'^{bin}| := n_2$, then their distance is the smaller possible one, *i.e.*:

$$h(\mathbf{X}_t^{bin}, \overline{\mathbf{x}}'^{bin}) = \min_{1 \leq i \leq n_1 - n_2 + 1} \left\{ h(\mathbf{X}_i, \overline{\mathbf{x}}'^{bin}) \,|\, \mathbf{X}_i := \left( X_{t,i}^{bin}, X_{t,i+1}^{bin}, \ldots, X_{t,i+n_2-1}^{bin} \right) \right\},$$

$h$ as in Definition 16.

Table 7: Average number of mutations needed to reach $\overline{\mathbf{x}}'^{bin}$, for $N = 7$ and starting from a Hamming distance 5. In $\overline{\mathbf{x}}'^{bin}$, only 2 amino-acids are hydrophobic, so by Definition 16, the optimal Hamming distance one can reach is 2. We compare a model in which no deletions nor insertions are allowed and a model in which 10% of all mutations are deletions or insertions. We denote by $\widehat{\tau_{\{\overline{\mathbf{x}}'^{bin}\}}}_n$ the average value obtained over $n$ simulations and by $\widehat{\sigma}_n$ its corresponding estimated standard deviation.

| % deletions/insertions | $|\overline{\mathbf{x}}'^{bin}|$ | $h(\mathbf{X}_0^{bin}, \overline{\mathbf{x}}'^{bin})$ | $n$ | $\widehat{\tau_{\{\overline{\mathbf{x}}'^{bin}\}}}_n$ | $\frac{\widehat{\sigma}_n}{\sqrt{n}}$ |
|---|---|---|---|---|---|
| **0** | 7 | 5 | 5000 | 374.28 | 5.38 |
| **10** | 7 | 5 | 5000 | 251.48 | 3.54 |

In this case, and thanks to the definition of Hamming distance as the minimal one, we clearly have more chances to obtain a good B-cell trait. This is confirmed by the results collected in Table 7. When deletions and insertions can occur, even with very weak probability, and if we allowed the BCR length to be greater than the antigen one, then the expected number of mutations needed to built the optimal BCR is more than 30% smaller.

## 5 Conclusion

In this paper, we have introduced a mathematical framework to study the impact of various mutation rules on the exploration of the space of traits in an evolutionary model. In particular, we have connected mutation rules to characteristic time-scales, such as hitting-times, through the study of associated graph structures. As a leading example, which was the original motivation for this study, we have considered applications of these results to the modeling of somatic hypermutations in the germinal center. The models considered so far do not include division and selection, which would lead to studying branching random walks on graphs, a topic of ongoing research.